\documentclass[12pt]{amsart}
\usepackage{amsfonts}
\usepackage{amssymb}
\usepackage{amsmath}
\input xy
\xyoption {all}

\usepackage{geometry}
\usepackage{amsthm}
\usepackage{amssymb}
\usepackage{latexsym}
\usepackage[dvips]{graphics}

\newtheorem{thm}{Theorem}[section] 
\newtheorem{cor}[thm]{Corollary}
\newtheorem{defn}[thm]{Definition}
\newtheorem{example}[thm]{Example}

\newtheorem{lemma}[thm]{Lemma}
\newtheorem{prop}[thm]{Proposition}
\newtheorem{remark}[thm]{Remark}

\newtheorem*{ack}{Acknowledgements}

\numberwithin{equation}{section}

\newcommand{\DC}{{\mathcal D}}
\newcommand{\LL}{{\mathcal L}}

\newcommand{\OO}{{\mathcal O}}
\newcommand{\HC}{{\mathcal H}}
\newcommand{\VV}{{\mathcal V}}

\newcommand{\Z}{\mathbb{Z}}

\newcommand{\Q}{\mathbb{Q}}
\newcommand{\R}{\mathbb{R}}
\newcommand{\C}{\mathbb{C}}

\newcommand{\HB}{\mathbb{H}}
\newcommand{\PP}{\mathbb{P}}

\newcommand{\VB}{\mathbb{V}}

\begin{document}
\date{\today}
\title{Hodge genera of algebraic varieties, II.}

\author[Sylvain E. Cappell ]{Sylvain E. Cappell}
\address{S. E. Cappell: Courant Institute, New York University, New York, NY-10012}
\email {cappell@cims.nyu.edu}

\author[Anatoly Libgober ]{Anatoly Libgober}
\address{A. Libgober : Department of Mathematics,
          University of Illinois at Chicago,
          851 S Morgan Street, Chicago, IL 60607,
                 USA.}
\email {libgober@math.uic.edu}

\author[Laurentiu Maxim ]{Laurentiu Maxim}
\address{L. Maxim : Department of Mathematics,
          University of Illinois at Chicago,
          851 S Morgan Street, Chicago, IL 60607,
                 USA.}
\email {lmaxim@math.uic.edu}

\author[Julius L. Shaneson ]{Julius L. Shaneson}
\address{J. L. Shaneson: Department of Mathematics, University of Pennsylvania, Philadelphia, PA-19104}
\email {shaneson@sas.upenn.edu}

\subjclass[2000]{ Primary 57R20, 32S20, 14C30, 32S35, 32S50, 14D05,
14D06, 14D07, 55N33; Secondary  57R45, 32S60, 13D15, 16E20. }

\thanks{A. Libgober partially supported by an NSF grant. S. Cappell and J. Shaneson partially supported by grants from DARPA}

\begin{abstract} We study the behavior of Hodge-theoretic genera
under morphisms of complex algebraic varieties. We prove that the
additive $\chi_y$-genus which arises in the motivic context
satisfies the so-called ``stratified multiplicative property", which
shows how to compute the invariant of the source of a proper
surjective morphism from its values on various varieties that arise
from the singularities of the map. By considering morphisms to a
curve, we obtain a Hodge-theoretic analogue of the Riemann-Hurwitz
formula.  We also study the contribution of monodromy to the
$\chi_y$-genus of a smooth projective family, and prove an
Atiyah-Meyer type formula for twisted $\chi_y$-genera. This formula
measures the deviation from multiplicativity of the $\chi_y$-genus,
and expresses the correction terms as higher-genera associated to
cohomology classes of the quotient of the total period domain by the
action of the monodromy group. By making use of Saito's theory of
mixed Hodge modules, we also obtain formulae of Atiyah-Meyer type
for the corresponding Hirzebruch characteristic classes.
\end{abstract}

\maketitle

\tableofcontents


\section{Introduction}

In the mid $1950$'s, Chern, Hirzebruch and Serre \cite{CHS} proved
that if $F \rightarrow E \rightarrow B$ is a fiber bundle of closed,
oriented, topological manifolds such that the fundamental group of
$B$ acts trivially on the cohomology of $F$, then the signatures of
the spaces involved are related by a simple multiplicative relation:
$\sigma(E)=\sigma(F)\sigma(B)$. A decade later, Atiyah \cite{At},
and respectively Kodaira \cite{Ko}, observed that without the
assumption on the action of the fundamental group the
multiplicativity relation fails. Moreover, Atiyah  showed that the
deviation from multiplicativity is controlled by the cohomology of
the fundamental group of $B$.

The main goal of this paper is to describe in a systematic way
multiplicativity properties of the Hirzebruch $\chi_y$-genus (and
the associated characteristic classes) for fibrations of algebraic
manifolds. We extend the results of Chern, Hirzebruch and Serre in
several different directions. First, we show that in the case when
the Chern-Hirzebruch-Serre assumption of the triviality of the
monodromy action is fulfilled, the $\chi_y$-genus is multiplicative.
Since by the Hodge index theorem, the signature of a K\"ahler
manifold is one of the values of the $\chi_y$-genus, this theorem
can be viewed as an extension of the Chern-Hirzebruch-Serre result
in the algebraic case. Secondly, we consider fibrations with
non-trivial monodromy action, and prove a Hodge-theoretic analogue
of the Atiyah signature formula. We also derive a formula for the
$\chi_y$-genus of $E$ in which the correction from the
multiplicativity of the $\chi_y$-genus is measured via pullbacks
under the period map associated with our fibration of certain
cohomology classes of the quotient of the period domain by the
action of the monodromy group. For some manifolds $F$ serving as a
fiber of the fibration in discussion,  as well as in the case when
one is interested in the value of polynomial $\chi_y$ yielding the
signature, the quotient of the period domain is the classifying
space of the monodromy group and our correction terms coincide with
those of Atiyah. In fact, the Atiyah terms are the appropriate
Novikov type higher-signatures, and our correction terms are
extensions of these Novikov invariants to the algebraic category. We
indicate that recently established birational properties of higher
genera \cite{BW,BL,R} are  valid also for ``the corrections to
multiplicativity" that we introduce here. Analogous formulae for the
associated Hirzebruch characteristic classes are also discussed.

\bigskip

We now present in more detail the content of each section and
summarize our main results.

In Section \ref{def}, we study the behavior of $\chi_y$-genera under
maps of complex algebraic varieties. We first consider a morphism
$f:E \to B$ of complex algebraic varieties with $B$ smooth and
connected, which is a fiber bundle in the (strong) complex topology,
and show that under certain assumptions on monodromy the
$\chi_y$-genera are multiplicative (cf. Lemma \ref{mult} and Lemma
\ref{cmult}). Such multiplicativity properties of genera were
previously studied in certain special cases in connection with
rigidity (e.g., see \cite{H,HBJ,Oc}). For instance, Hirzebruch's
$\chi_y$-genus is multiplicative in bundles of (stably) almost
complex manifolds with structure group a compact connected Lie group
(the latter condition implies trivial monodromy), and in fact it is
uniquely characterized by this property. The proof of our
multiplicativity result uses the fact that the Leray spectral
sequences of the map $f$ are  spectral sequences in the category of
mixed Hodge structures. The latter claim is a consequence of Saito's
theory of mixed Hodge modules (e.g., see \cite{Sa1}), and is
discussed in some detail in \S \ref{speq}.

In Section \ref{Sing}, we consider proper morphisms that are allowed
to have singularities, and extend the above multiplicativity
property to this general stratified case. More precisely, we prove
that, under the assumption of trivial monodromy along the strata of
our map, the additive $\chi_y$-genus that arises in the motivic
context satisfies the so-called ``stratified multiplicative
property" (cf. Proposition \ref{e} and Corollary \ref{ee}). This
property shows how to compute the invariant of the source of a
proper surjective morphism from its values on various varieties that
arise from the singularities of the map, thus yielding powerful
topological constraints on the singularities of any algebraic map.
It also provides a method of inductively computing these genera of
varieties. A similar result was obtained by Cappell, Maxim and
Shaneson, for the behavior of intersection homology Hodge-theoretic
invariants, both genera and characteristic classes (see \cite{CMS},
and also \cite{CMS0}). Such formulae were first predicted by Cappell
and Shaneson in the early $1990$'s, see the announcements \cite{CS,
S}, following their earlier work on stratified multiplicative
properties for signatures and associated topological characteristic
classes defined using intersection homology, cf. \cite{CS2} (see
also [\cite{BSY}, \S 4], and \cite{Y2} for a functorial
interpretation of Cappell-Shaneson's $L$-classes).

In the special case of maps to a smooth curve, and under certain
assumptions for the monodromy along the strata of special fibers, in
\S \ref{RH} we obtain a Hodge-theoretic analogue of the
Riemann-Hurwitz formula \cite{K}. The proof uses Hodge-theoretic
aspects of the nearby and vanishing cycles in the context of
one-parameter degenerations of projective varieties.

The contribution of monodromy to $\chi_y$-genera is studied in
Section \ref{nm}. This can be applied to compute the summands
arising from singularities in the above formulae. For simplicity, we
first consider a smooth proper map $f:E \to B$ of smooth projective
varieties (thus a fibration in the strong topology), and compute
$\chi_y(E)$ so that the (monodromy) action of $\pi_1(B)$ on the
cohomology of the typical fiber is taken into account (see Theorem
\ref{AM}). The proof uses the Hirzebruch-Riemann-Roch theorem
\cite{H} and standard facts from the classical Hodge theory. Our
formula (\ref{twisted}) is a Hodge-theoretic analogue of Atiyah's
formula for the signature of fiber bundles \cite{At}, and measures
the deviation from multiplicativity of the $\chi_y$-genus in the
presence of monodromy. In Section \ref{periods}, this deviation is
expressed in terms of higher-genera associated to cohomology classes
of the quotient of the total period domain by the action of the
monodromy group.  As a corollary of our formula (\ref{twisted}), we
point out that if the action of $\pi_1(B)$ on the cohomology of the
typical fiber $F$ preserves the Hodge filtration, then the
$\chi_y$-genus is still multiplicative, i.e., $\chi_y(E)=\chi_y(B)
\cdot \chi_y(F)$. This assertion is false in the non-compact case
(see Example \ref{counterex} (2)). As a byproduct of the proof of
Theorem \ref{AM}, we obtain a Hodge-theoretic analogue of Meyer's
formula for twisted signatures \cite{Me} (see Corollary
\ref{Meyer}). At the end of Section \ref{nm}, we present several
interesting extensions of these Hodge-theoretic Atiyah-Meyer
formulae to more general situations, where we allow the spaces
involved to be singular and/or non-compact.

In Section \ref{CLASS}, we extend some of the above mentioned
results on $\chi_y$-genera to Atiyah-Meyer type formulae for the
corresponding Hirzebruch characteristic classes. The proofs are much
more involved, and use in an essential way the construction of
Hirzebruch classes via Saito's theory of mixed Hodge modules (cf.
\cite{BSY}). At this end, we point out that many of our formulae for
genera can be obtained from the corresponding formulae for
characteristic classes. However, we chose to present the results for
genera first, since they can be proven by using only standard facts
of classical Hodge theory, thus becoming accessible to a wider
audience.

It is conceivable that many of our results remain valid in a more
general context (e.g., for compactifiable complex analytic
varieties). However, for simplicity, we present here our
Atiyah-Meyer type results in the algebraic setting.

We have tried to make this paper as self-contained as possible. For
this reason, in \S \ref{Saito} we provide necessary background on
Saito's theory of mixed Hodge modules, and in \S \ref{Van} we recall
Deligne's formalism of nearby and vanishing cycles. However, we
assume reader's familiarity with certain aspects of Deligne's Hodge
theory (\cite{De, PS}).

In a future paper, we will consider extensions of our monodromy
formulae to the singular setting, both for genera and characteristic
classes. One possible such extension makes use of intersection
homology \cite{GM1,GM2} and the BBDG decomposition theorem
\cite{BBD,Cat}. This approach is motivated by the considerations in
\cite{CMS} (where the case of trivial monodromy was considered), and
by an extension of the Atiyah-Meyer signature formula to the
singular case, which is due to Banagl, Cappell and Shaneson
\cite{BCS}.

\begin{ack} We are grateful to J\"org Sch\"urmann for reading a
first draft of this work, and for making many valuable comments and
suggestions for improvement. We thank Mark Andrea de Cataldo,
Alexandru Dimca and Nero Budur for many inspiring conversations on
this subject, and Joseph Steenbrink and Shoji Yokura for commenting
on a preliminary version of this work.
\end{ack}

\section{Hodge genera and singularities of maps}\label{def}

\subsection{Hodge genera. Definitions.}\label{Defns} In this section,
we define the Hodge-theoretic invariants of complex algebraic varieties,
which will be studied in the sequel. We assume reader's familiarity
with Deligne's theory of mixed Hodge structures \cite{De}.

For any complex algebraic variety $Z$, we define the
$\chi_y^c$-genus in terms of the Hodge-Deligne numbers of compactly
supported cohomology of $Z$ (cf. \cite{D}). More precisely,
$$\chi_y^c(Z)=\sum_p \left( \sum_{i,q} (-1)^{i-p} h^{p,q}(H^i_c(Z;\C)) \right) y^p =
\sum_{i,p \geq 0} (-1)^{i-p} \text{dim}_{\C} Gr^p_F H^i_c(Z,\C)
\cdot y^p,$$ where $h^{p,q}(H^i_c(Z;\C))=\text{dim}_{\C} Gr^p_F
(Gr^W_{p+q} H^i_c(Z) \otimes \C)$, with $F^{\bullet}$ and
$W_{\bullet}$ the Hodge and respectively the weight filtration of
Deligne's mixed Hodge structure on $H^i_c(Z)$. Similarly, we define
the $\chi_y$-genus of $Z$, $\chi_y(Z)$, by using the Hodge-Deligne
numbers of $H^*(Z;\C)$. Of course, for a complete variety $Z$ we
have that $\chi_y(Z)=\chi_y^c(Z)$. If $Z$ is smooth and projective,
then each cohomology group $H_c^i(Z;\C)=H^i(Z;\C)$ has a pure Hodge
structure of weight $i$, and the above formulae define Hirzebruch's
$\chi_y$-genus (cf. \cite{H}). Note that for any complex variety
$Z$, we have that $\chi^c_{-1}(Z)=\chi_{-1}(Z)=\chi(Z)$ is the usual
Euler characteristic, where for the first equality we refer to
\cite{F}, pp. 141-142. Similarly, $\chi_0$ and $\chi_0^c$ are two
possible extensions to singular varieties of the arithmetic genus.

The compactly supported $\chi_y$-genus, $\chi_y^c$, satisfies the so-called ``scissor relations" for complex
varieties, that is: $\chi_y^c(Z)=\chi_y^c(W)+\chi_y^c(Z \setminus
W)$, for $W$ a closed subvariety of $Z$. Therefore, $\chi_y^c$ can
be defined on $K_0(Var_{\C})$, the Grothendieck group of varieties
over $\C$ which arises in the motivic context.

More generally, we can define $\chi_y$-genera on the Grothendieck
group of mixed Hodge structures $K_0(mhs)=K_0(D^b mhs)$, where we
denote by $mhs$ the abelian category of (rational) mixed Hodge
structures. Indeed, if $K \in mhs$, define
\begin{equation}\label{ab} \chi_y([K]):=\sum_p \text{dim}_{\C} Gr^p_F (K \otimes
\C) \cdot (-y)^p,\end{equation} where $[K]$ is the class of $K$ in
$K_0(mhs)$. This is well-defined on  $K_0(mhs)$  since the functor
$Gr^p_F$ preserves exactness. For $K^{\bullet}$  a bounded complex
of mixed Hodge structures, we define
$$[K^{\bullet}]:=\sum_{i \in \Z} (-1)^i [K^i] \ \ \in K_0(mhs)$$
and note that we have: $$[K^{\bullet}]=\sum_{i \in \Z} (-1)^i
[H^i(K^{\bullet})].$$ In view of (\ref{ab}), we set
\begin{equation}\label{defn}\chi_y([K^{\bullet}]):=\sum_{i \in \Z} (-1)^i
\chi_y([K^i]).\end{equation} In this language, we have that:
$$\chi_y^c(Z)=\chi_y([H_c^{\bullet}(Z;\Q)])$$
and
$$\chi_y(Z)=\chi_y([H^{\bullet}(Z;\Q)]),$$
where $H_c^{\bullet}(Z;\Q)$ and $H^{\bullet}(Z;\Q)$ are regarded as
bounded complexes of mixed Hodge structures, with all differentials
equal to zero.

\subsection{Basics of Saito's theory of mixed Hodge modules.}\label{Saito} Even
though the theory of mixed Hodge modules is very involved, in this
section we give a brief overview adapted to our needs, and we show
some quick applications in the following sections. The standard
references for the algebraic case are Saito's papers \cite{Sa2} and
\cite{Sa1}, \S4, but see also the book \cite{PS}, Chapter $14$, for
a brief survey.

We recall that for any complex algebraic variety $Z$, the derived
category of bounded cohomologically constructible complexes of
sheaves of $\Q$-vector spaces on $Z$ is denoted by $D_c^b(Z)$, and
it contains as a full subcategory the category $\text{Perv}_{\Q}(Z)$
of perverse $\Q$-complexes. The Verdier duality operator
$\mathbb{D}_Z$ is an involution on $D_c^b(Z)$ preserving
$\text{Perv}_{\Q}(Z)$. Associated to a morphism $f:X \to Y$ of
complex algebraic varieties, there are pairs of adjoint functors
$(f^*, Rf_*)$ and $(f^!, Rf_!)$ between the respective categories of
cohomologically constructible complexes, which are interchanged by
Verdier duality. For details, see the books \cite{Di2, Sc}.

M. Saito associated to a complex algebraic variety $Z$ an abelian
category $MHM(Z)$, the category of \emph{mixed Hodge modules} on
$Z$, together with a faithful forgetful functor
$$rat: D^bMHM(Z) \to D^b_c(Z)$$
such that $rat(MHM(Z)) \subset \text{Perv}_{\Q}(Z)$. For
$M^{\bullet} \in D^bMHM(Z)$, $rat(M^{\bullet})$ is called the
underlying rational complex of $M^{\bullet}$.

We say that $M \in MHM(Z)$ is supported on $S$ if and only if
$rat(M)$ is supported on $S$.  Saito showed that the category of
mixed Hodge modules supported on a point coincides with the category
of (graded) polarizable rational mixed Hodge structures. In this
case, the functor $rat$ associates to a mixed Hodge structure the
underlying rational vector space.

Since $MHM(Z)$ is an abelian category, the cohomology groups of any
complex $M^{\bullet} \in D^bMHM(Z)$ are mixed Hodge modules. The
underlying rational complexes of the cohomology groups of a complex
of mixed Hodge modules are the perverse cohomologies of the
underlying rational complex, that is,
$rat(H^j(M^{\bullet}))={^p\HC}^j(rat(M^{\bullet}))$.

The Verdier duality functor $\mathbb{D}_Z$ lifts to $MHM(Z)$ as an
involution, in the sense that it commutes with the forgetful
functor: $rat \circ \mathbb{D}_Z = \mathbb{D}_Z \circ rat$.

For a morphism $f:X \to Y$ of complex algebraic varieties, there are
induced functors $f_*, f_! : D^bMHM(X) \to D^bMHM(Y)$ and $f^*, f^!
: D^bMHM(Y) \to D^bMHM(X)$, exchanged under the Verdier duality
functor, and which lift the analogous functors on the level of
constructible complexes. Moreover, if $f$ is proper, then $f_!=f_*$.

Let us give a rough picture of what the objects in Saito's category
of mixed Hodge modules look like. For $Z$  \emph{smooth}, $MHM(Z)$
is a full subcategory of the category of objects $((M,F), K, W)$
such that:
\begin{enumerate} \item $(M,F)$ is an algebraic holonomic filtered $\DC$-module $M$ on
$Z$, with an increasing ``Hodge" filtration $F$ by coherent
algebraic $\OO_Z$-modules; \item $K \in \text{Perv}_{\Q}(Z)$ is the
underlying rational sheaf complex, and there is a quasi-isomorphism
$\alpha: DR(M) \simeq \C \otimes K$ in $\text{Perv}_{\C}(Z)$, where
$DR$ is the de Rham functor shifted by the dimension of $Z$; \item
$W$ is a pair of (weight) filtrations on $M$ and $K$ compatible with
$\alpha$.
\end{enumerate} For a singular variety $Z$, one works with local embeddings
into manifolds and corresponding filtered $\DC$-modules with support
on $Z$. In this notation, the functor $rat$ is defined by
$rat((M,F), K, W)=K$.

A complex  $M^{\bullet} \in D^bMHM(Z)$ is \emph{mixed of weight
$\leq k$ (resp. $\geq k$)} if $Gr_i^W H^jM^{\bullet} = 0$ for all $i
> j+k$ (resp. $i < j+k$), and it is \emph{pure of weight $k$} if
$Gr_i^W H^jM^{\bullet} = 0$ for all $i \neq j+k$. If $f$ is a map of
algebraic varieties, then $f_!$ and $f^*$ preserve weight $\leq k$,
and $f_*$ and $f^!$ preserve weight $\geq k$. If $M^{\bullet} \in
D^bMHM(Z)$ is of weight $\leq k$ (resp. $\geq k$), then
$H^jM^{\bullet}$ has weight $\leq j+k$ (resp. $\geq j+k$).

If $j: U \hookrightarrow Z$ is a Zariski-open subset in $Z$, then
the intermediate extension $j_{!*}$ (cf. \cite{BBD}) preserves the
weights.

Following \cite{Sa1}, there exists a unique object $\Q^H \in
MHM(point)$ such that $rat(\Q^H)=\Q$ and $\Q^H$ is of type $(0,0)$.
In fact, $\Q^H=((\C,F), \Q, W)$, with $gr^F_i=0=gr^W_i$ for all $i
\neq 0$, and $\alpha : \C \to \C \otimes \Q$ the obvious
isomorphism. For a complex variety $Z$, define $\Q_Z^H:=a_Z^*\Q^H
\in D^bMHM(Z)$ with $a_Z:Z\to point$ the map to a point. If $Z$ is
\emph{smooth} of dimension $n$, then $\Q_Z[n] \in Perv(\Q_X)$ and
$\Q_Z^H[n]\in MHM(Z)$ is a single mixed Hodge module (in degree
$0$), explicitely described by  $\Q_Z^H[n]=((\OO_Z, F), \Q_Z[n],
W),$ where $F$ and $W$ are trivial filtrations so that
$gr^F_i=0=gr^W_{i+n}$ for all $ i \neq 0$. So if $Z$ is smooth of
dimension $n$, then $\Q_Z^H[n]$ is pure of weight $n$. By the
stability of the intermediate extension functor, this shows that if
$Z$ is any algebraic variety and $j:U \hookrightarrow Z$ is the
inclusion of a smooth Zariski-open subset, then the intersection
cohomology module $IC_Z^H:=j_{!*}(\Q_U^H[n])$ is pure of weight $n$.

More generally, if $\mathbb{V}$ is a polarized variation of Hodge
structures of weight $k$ with quasi-unipotent monodromy at
infinity\footnote{This condition is automatically satisfied if $\VB$
is a variation of $\Z$-Hodge structures, by the monodromy theorem,
\cite{PS}, Thm. 11.8.} defined on a Zariski-open subset $U$ of $Z$,
then $\mathbb{V}$ corresponds to a smooth mixed Hodge module
$\mathbb{V}^H[n]$ on $U$ (i.e., the associated rational complex  is
a local system) of weight $k+n$, whose underlying perverse sheaf is
$\VB[n]$ (see \cite{Sa0}, Thm. 5.4.3, \cite{Sa2}, \S2, \cite{PS},
Thm. 14.30). So the twisted (middle-perversity) intersection
homology complex $IC_Z^H(\mathbb{V}):=j_{!*}(\VB^H[n])$ is a mixed
Hodge module, pure of weight $k+n$.

The following result is implicit in the work of Saito.

\begin{prop}\label{Sai} Let $Z$ be a $n$-dimensional irreducible projective
variety, and $\VB$ a polarized variation of Hodge structures of
weight $k$ (with quasi-unipotent monodromy at infinity) defined on a
Zariski-open dense subset of $Z$. Then the intersection cohomology
group $IH^j(Z; \VB)$ carries a pure Hodge structure of weight $j+k$.
\end{prop}

\begin{proof} Let $a_Z:Z \to point$ be the constant map to a point.
Since $a_Z$ is proper, ${a_Z}_*$ preserves the weights. The claim
follows from the isomorphism
$$IH^j(Z;\VB) \cong \HB^{j-n}(Z;IC_Z^H(\mathbb{V})) \cong
H^{j-n}({a_Z}_* IC_Z^H(\mathbb{V})),$$ by noting that ${a_Z}_*
IC_Z^H(\mathbb{V})$ is a pure complex of weight $n+k$, thus
$H^{j-n}({a_Z}_* IC_Z^H(\mathbb{V}))$ is a pure Hodge module of
weight $j+k$ supported over a point. Since over a point mixed Hodge
modules are exactly the (graded) polarizable mixed Hodge structures,
the latter is a Hodge structure of weight $j+k$.

\end{proof}

The following corollary will be needed in the sequel:
\begin{cor}\label{Hp} Let $Z$ be a smooth complex projective
algebraic variety, and $\VB$ a polarized variation of Hodge
structures of weight $k$ defined on $Z$. Then the cohomology group
$H^j(Z; \VB)$ carries a pure Hodge structure of weight $j+k$.
\end{cor}

\begin{remark}\label{mixed}\rm More generally,
if the variety $Z$ in Proposition \ref{Sai} and Corollary \ref{Hp}
is not necessarily compact, and if $\VB$ is a polarized variation of
Hodge structures or, more generally, an admissible variation of
mixed Hodge structures (for a definition, see \cite{PS}, Def. 14.47
and the references therein) on (a Zariski-open dense subset of) $Z$,
then the associated (intersection) cohomology groups carry natural
mixed Hodge structures.
\end{remark}

\subsection{Spectral sequences of mixed Hodge modules.}\label{speq}

In this section, we justify the claim that certain Leray-type
spectral sequences (e.g., the Leray spectral sequence of an
algebraic morphism, or the hypercohomology spectral sequence) are in
fact spectral sequences of mixed Hodge structures.

From the general theory of spectral sequences, since the category of
mixed Hodge modules is abelian, the canonical filtration $\tau$ on
$D^bMHM(Z)$ preserves complexes of mixed Hodge modules. Therefore,
the second fundamental spectral sequence (\cite{PS}, \S A.3.4) for
any (left exact) functor $F$ sending mixed Hodge modules to mixed
Hodge modules, that is, the spectral sequence
\begin{equation}\label{gen}
E_2^{p,q}=H^pF(H^q(M^{\bullet})) \Longrightarrow
H^{p+q}F(M^{\bullet}),
\end{equation}
is a spectral sequences of mixed Hodge modules.

Note that the canonical $t$-structure $\tau$ on $D^bMHM(Z)$
corresponds to the perverse truncation $^p\tau$ on $D^b_c(Z)$.
However, Saito (\cite{Sa1}, Remark 4.6(2)) constructed another
$t$-structure $'\tau$ on $D^bMHM(Z)$ that corresponds to the
classical $t$-structure on $D^b_c(Z)$. By using the $t$-structure
$'\tau$ in the construction of the second fundamental spectral
sequence above, one can show that the classical Leray spectral
sequences are in fact spectral sequences of mixed Hodge structures.

\begin{example}\rm\label{LSS} \emph{Leray spectral sequences.}\newline
$(1)$ \ Let $Z$ be a complex algebraic variety. Then for
$\mathcal{F}^{\bullet}$ a bounded complex of sheaves with
constructible cohomology on $Z$, we have the spectral sequence with
the $E_2$-term given by
\begin{equation}\label{leray}
E_2^{p,q}=H^p(Z;\mathcal{H}^q(\mathcal{F}^{\bullet}))
\Longrightarrow \mathbb{H}^{p+q}(Z;\mathcal{F}^{\bullet}).
\end{equation}
This spectral sequence is induced by the natural filtration on the
complex $\mathcal{F}^{\bullet}$, and if $\mathcal{F}^{\bullet}$
underlies a complex of mixed Hodge modules, then the spectral
sequence is compatible with mixed Hodge structures: this follows by
using the $t$-structure $'\tau$ defined in \cite{Sa1}, Remark
4.6(2), and the fundamental spectral sequence (\ref{gen}) for
$F=\Gamma(Z,\cdot)=({a_Z})_*$, together with the fact that mixed
Hodge modules over a point are (graded polarizable) mixed Hodge
structures.

Similarly, by taking $F=\Gamma_c(Z,\cdot)=({a_Z})_!$ together with
the $t$-structure $'\tau$ above, the compactly supported
hypercohomology Leray spectral sequence
\begin{equation}\label{csleray}
E_2^{p,q}=H_c^p(Z;\mathcal{H}^q(\mathcal{F}^{\bullet}))
\Longrightarrow \mathbb{H}_c^{p+q}(Z;\mathcal{F}^{\bullet}).
\end{equation}
is a spectral sequence in the category of mixed Hodge structures,
provided $\mathcal{F}^{\bullet}$ underlies a bounded complex of
mixed Hodge modules.\newline

\noindent $(2)$ \ Let $f:E \to B$ be a morphism of complex algebraic
varieties. The Leray spectral sequence for $f$
\begin{equation}\label{mapleray}E_2^{p,q}=H^p(B, R^qf_*\Q_E)
\Longrightarrow H^{p+q}(E)\end{equation} is a spectral sequence of
mixed Hodge structures: this follows by using $'\tau$ on
$D^bMHM(B)$, and the example above applied to
$\mathcal{F}^{\bullet}=Rf_*\Q_E$, that underlies
$M^{\bullet}=f_*\Q^H_E \in D^bMHM(B)$. Similarly, there is a
compactly supported version of the Leray spectral sequence for $f$
in the category of mixed Hodge structures, namely
\begin{equation}\label{cmapleray}E_2^{p,q}=H_c^p(B, R^qf_!\Q_E)
\Longrightarrow H_c^{p+q}(E).\end{equation}
By using the usual $t$-structure on $D^bMHM(B)$, hence the perverse $t$-structure on $D^b_c(B)$, we obtain that the perverse Leray spectral sequence for $f$,
\begin{equation}\label{pL}
E^{i,j}_2=\HB^i(B, {^p\HC}^j( f_*\mathcal{F}^{\bullet} )) \Longrightarrow \HB^{i+j}(E; \mathcal{F}^{\bullet})
\end{equation}
is a spectral sequence in the category of mixed Hodge structures,
provided $\mathcal{F}^{\bullet}$ underlies a bounded complex of
mixed Hodge modules.

\end{example}

\begin{remark}\rm In the quasi-projective setting, a different proof
of the fact that the spectral sequences (\ref{mapleray}) and
(\ref{cmapleray}) are spectral sequences of mixed Hodge structures
was given by Arapura in \cite{Ara}.
\end{remark}

\subsection{Multiplicativity properties of $\chi_y$-genera.}\label{lem} In this
section we use the above spectral sequences in studying
multiplicative properties of the Hodge-theoretic genera. We first
prove the following:
\begin{lemma}\label{mult} Let $E$, $B$, $F$ be complex algebraic varieties with $B$ smooth and connected, and
let $p:E \to B$ be an algebraic morphism such that $p$ is locally
trivial in the strong (complex) topology of $B$, with fiber $F$.
Assume that the local systems $R^jp_*\C_E$ are constant for each $j$
(e.g., $\pi_1(B)=0$). Then
$$\chi_y(E)=\chi_y(B)\chi_y(F).$$
If moreover, all $E$, $B$ and $F$ are smooth, then
$$\chi_y^c(E)=\chi_y^c(B)\chi_y^c(F).$$
\end{lemma}
\begin{proof}
Consider the Leray spectral sequence of $p$, that is
$$E_2^{p,q}=H^p(B, R^qp_*\C_E) \Longrightarrow H^{p+q}(E)$$
In the category of algebraic varieties this is a spectral sequence
of mixed Hodge structures (cf. \S \ref{speq}).

Since the local systems $R^qp_*\C_E$ are constant, the corresponding
variations of mixed Hodge structures are trivial as the natural
morphism $H^0(B, R^qp_*\C_E) \to (R^qp_*\C_E)_y$ is a mixed Hodge
structure isomorphism for any $y \in B$. Then we have $H^p(B,
R^qp_*\C_E)=H^p(B) \otimes (R^qp_*\C_E)_y$ as mixed Hodge
structures, and there are mixed Hodge structure isomorphisms
\begin{equation}\label{Lt}E_2^{p,q}=H^p(B) \otimes H^q(F)\end{equation}

Since all differentials in the Leray spectral sequence are mixed
Hodge structure morphisms, thus strict with respect to the Hodge and
weight filtrations, by \cite{De}, Lemme 1.1.11 we get a spectral
sequence for the Hodge components of a given type $(k,l)$:
\begin{equation}\label{Lh} E(k,l)_2^{p,q}:=Gr^k_F Gr^W_{k+l}E_2^{p,q}
\Longrightarrow Gr^k_F Gr^W_{k+l}H^{p+q}(E)
\end{equation}
Now let $e^{k,l}$ be the Euler characteristic of Hodge-type $(k,l)$,
i.e., for a  complex algebraic variety $Z$ we define
$$e^{k,l}(Z)=\sum_i(-1)^ih^{k,l}(H^i(Z)).$$
By the invariance of Euler characteristics under spectral sequences,
from (\ref{Lt}) and (\ref{Lh}) we obtain
\begin{eqnarray*} e^{k,l}(E) &=& \sum_{i}
(-1)^i\text{dim}\left(\oplus_{p+q=i} E(k,l)_2^{p,q}\right) \\
&=& \sum_{r+t=k, s+u=l} e^{r,s}(B) e^{t,u}(F). \end{eqnarray*} The
multiplicativity of $\chi_y$ follows by noting that for a variety
$Z$,
$$\chi_y(Z)=\sum_{k,l}e^{k,l}(Z) \cdot (-y)^k.$$
The claim about the multiplicativity of the $\chi_y^c$-genus follows
by Poincar\'e Duality, by noting that if $Z$ is a smooth complex
algebraic variety of dimension $n$, then
$\chi_y^c(Z)=(-y)^n\chi_{y^{-1}}(Z)$. Indeed, the Poincar\'e duality
isomorphism takes classes of type $(p,q)$ in $H^j_c(Z)$ to classes
of type $(n-p,n-q)$ in $H^{2n-j}(Z)$, where $n=\text{dim} Z$.

\end{proof}

By noting that for any $b \in B$, $(R^jp_!\C_E)_b \cong
H^j_c(F;\C)$, a similar argument applied to the compactly supported
Leray spectral sequence (\ref{cmapleray}) of the map $p$, yields the
following

\begin{lemma}\label{cmult} Let $E$, $B$, $F$ be complex algebraic varieties with $B$ smooth and connected,
and consider $p:E \to B$ an algebraic morphism such that $p$ is
locally trivial in the strong (complex) topology of $B$ with fibre
$F$. Assume that the local systems $R^jp_!\C_E$ are constant for
each $j$ (e.g., $\pi_1(B)=0$). Then
$$\chi^c_y(E)=\chi^c_y(B)\chi^c_y(F).$$
\end{lemma}

\begin{remark}\label{E}\rm
The same argument can be used to show that the results of the above
lemmas hold for the Hodge-Deligne polynomials (or the $E$-functions)
defined by $E(Z;u,v)=\sum_{k,l} e^{k,l}(Z)u^kv^l$ (and similarly for
the $E$-functions $E_c(Z;u,v)$ defined by using the compactly
supported cohomology).\footnote{Note that $\chi_y(Z)=E(Z;-y,1)$, and
similarly, $\chi_y^c(Z)=E_c(Z;-y,1)$.} In particular, these results
hold for the weight polynomials $W(Z;t):=E(Z;t,t)$ and respectively
$W_c(Z;t):=E_c(Z;t,t)$ considered in \cite{DL}. In fact, Lemma
\ref{mult} is modeled after \cite{DL}, Theorem 6.1.
\end{remark}

\begin{example}\rm\label{counterex}
$(1)$ \ As an example, consider the case of the Hopf fibration
defining $\C\PP^n$. Then
$\chi_y(\C\PP^n)=\chi_y^c(\C\PP^n)=1+(-y)+\cdots +(-y)^n$,
$\chi_y^c(\C^{n+1} \setminus \{0\})=(-y)^{n+1}-1$,
$\chi_y^c(\C^*)=-y-1$, and by Poincar\'e Duality, $\chi_y(\C^{n+1}
\setminus \{0\})=1-(-y)^{n+1}$, $\chi_y(\C^*)=1+y$. Thus the
multiplicativity for both $\chi_y$ and $\chi_y^c$ holds.
\newline $(2)$ \ Let $p$ be the Milnor fibration of
a weighted homogeneous isolated hypersurface singularity at the
origin in $\C^{n+1}$, that is, $F=\{p=1\} \hookrightarrow E=\C^{n+1}
\setminus \{p=0\} \to B=\C^*$, for $p$ a weighted homogeneous
polynomial in $n+1$ variables, with an isolated singular point at
the origin. In this case, the monodromy is an algebraic morphism of
finite order, so it preserves the filtrations in $H^*(F;\C)$. In
other words, each $R^jp_*\Q_E$ is a (non-constant) local system of
mixed Hodge structures (note that this is not the case for general
singularities). The mixed Hodge structure on $H^*(F;\Q)$ is known by
work of Steenbrink \cite{St}. It turns out that even in this special
case, the $\chi_y$-genera are \emph{not} multiplicative (but see
Section \S \ref{nm} for the compact case). Here is a concrete
example: let $p(x,y)=x^3-y^2$ defining the cuspidal cubic in $\C^2$.
Then, in the notations above and by \cite{St}, the (mixed) Hodge
numbers of $F$ are $h^{0,0}(H^0(F))=1$, $h^{1,0}(H^1(F))=1$,
$h^{0,1}(H^1(F))=1$ and $h^{1,1}(H^1(F))=0$ (note that $H^2(F)=0$
since $F$ is affine of complex dimension $1$). Therefore, we obtain
that $\chi_y(F)=y$, so by Poincar\'e Duality it follows that
$\chi_y^c(F)=(-y)\cdot \chi_{y^{-1}}(F)=-1$. It also follows easily
that $\chi_y^c(E)=y^2+y$, $\chi_y^c(B)=-y-1$, and $\chi_y(E)=1+y$,
$\chi_y(B)=1+y$.

\end{example}

\begin{remark}\rm\label{zar}
The assumption of trivial monodromy is closely related to, but
different from the situation of ``algebraic piecewise trivial" maps
coming up in the motivic context (e.g., see \cite{BSY}). For
example, the formula in Lemma \ref{cmult} is true (without any
assumption on monodromy) for a Zariski locally trivial fibration of
possibly singular complex algebraic varieties (see \cite{D},
Corollary 1.9, or \cite{BSY}, Example 3.3). Note also that if all
spaces involved are smooth, then a Zariski locally trivial fibration
is a locally trivial fibration in the complex (strong) topology, and
the monodromy action is trivial.
\end{remark}

\subsection{$\chi^c_y$-genera and singularities of maps.}\label{Sing} In this
section, by analogy with the results of \cite{CMS0, CMS}, we discuss
the behavior of the $\chi^c_y$-genus under proper morphisms of
algebraic varieties, and show that $\chi_y^c$ satisfies the
\emph{stratified multiplicative property} in the sense of \cite{CS}.
The result will be further refined in \S \ref{curves}, in the case
of maps onto curves.

Let $f :X^n \to Y^m$ be a proper map of complex algebraic varieties
of indicated dimensions. Such a map can be stratified with
subvarieties as strata, i.e., there exist finite algebraic Whitney
stratifications $\mathcal{X}$ of $X$ and $\mathcal{S}$ of $Y$, such
that for any component $S$ of a stratum of $Y$, $f^{-1}(S)$ is a
union of connected components of strata of $\mathcal{X}$, each of
which is mapping submersively to $S$. This implies that
$f_{|f^{-1}(S)}: f^{-1}(S) \to S$ is a locally trivial map of
Whitney stratified spaces (see \cite{GM}, \S I.1.6). For simplicity,
we will assume that $f$ is smooth over the dense open stratum (e.g.,
the latter is connected). We denote by $F$ the general fiber of $f$,
and by $F_S$ the fiber of $f$ above the singular stratum $S \in
\mathcal{S}$.

We can now prove the following formula, showing the deviation from
multiplicativity of the $\chi_y^c$-genus in the case of a stratified
map (compare with results in \cite{CMS0, CMS}):

\begin{prop}\label{e} Let $f:X \to Y$ be a proper surjective morphism of (possibly
singular) complex algebraic varieties. Let $\mathcal{S}$ be the set
of components of open strata of $Y$ in a stratification of $f$, and
assume $\pi_1(S)=0$ for all $S \in \mathcal{S}$. For $S \in
\mathcal{S}$, define $\hat{\chi_y^c}(\bar{S})$ inductively by the
formula:
$$\hat{\chi_y^c}(\bar{S})=\chi_y^c(\bar S) - \sum_{W<S} \hat{\chi_y^c}(\bar{W}),$$
where the sum is over all $W \in \mathcal{S}$ with $W \subset \bar S
\setminus S$. Then:
\begin{equation}\label{chic}
\chi_y^c(X)=\chi^c_y(Y) \chi_y^c(F) + \sum_{S \in \mathcal{S},
\text{dim} S< \text{dim} Y} \hat{\chi_y^c}(\bar{S}) \cdot \left[
\chi^c_y(F_S)-\chi^c_y(F) \right].
\end{equation}

\end{prop}

\begin{proof}
Note that by the additivity of $\chi_y^c$-genera, we have that
$\hat{\chi_y^c}(\bar{S})=\chi^c_y(S)$. Next, by additivity and
multiplicativity for locally trivial fibrations with trivial
monodromy (cf. Lemma \ref{cmult}), it is easy to see that:
$$\chi^c_y(X)=\chi^c_y(Y) \chi^c_y(F) + \sum_{S \in \mathcal{S},  \text{dim} S< \text{dim} Y} \chi^c_y(S) \cdot
\left[ \chi^c_y(F_S)-\chi^c_y(F) \right].$$ The result follows.

\end{proof}

\begin{example}\rm \emph{Smooth blow-up}\newline
Let $Y$ be a smooth subvariety of codimension $r+1$ in a smooth
variety $X$.  Let $\pi:\tilde{X} \to X$ be the blow-up of $X$ along
$Y$.  Then $\pi$ is an isomorphism over $X \setminus Y$ and a
projective bundle (Zariski locally trivial) over $Y$, coresponding
to the normal bundle of $Y$ in $X$ of rank $r+1$. The formula in
Proposition \ref{e} yields
\begin{equation}\label{BSY}
\chi^c_y(\tilde{X})=\chi^c_y(X)+\chi^c_y(Y) \cdot \left( -y+\cdots
+(-y)^r \right).\end{equation} Note that this formula also holds
without any assumption on monodromy, by using instead Remark
\ref{zar} and the fact that $\pi^{-1}(Y)$ is a Zariski locally
trivial fibration over $Y$ with fiber $\C\PP^r$ (see \cite{D}, \S
1.10).
\end{example}

In the notations of Proposition \ref{e}, we also obtain the
following fact which was claimed in \cite{CMS}:

\begin{cor}\label{ee} Let $f:X \to Y$ be a proper surjective morphism of complex projective
algebraic varieties, and assume $\pi_1(S)=0$ for all $S \in
\mathcal{S}$. Then with the obvious definition for
$\hat{\chi_y}(\bar{S})$, the following holds:
\begin{equation}\label{chismp}
\chi_y(X)=\chi_y(Y) \chi_y(F) + \sum_{S \in \mathcal{S},  \text{dim}
S< \text{dim} Y} \hat{\chi_y}(\bar{S}) \cdot \left[
\chi_y(F_S)-\chi_y(F) \right].
\end{equation}
\end{cor}

\begin{remark}\rm More generally, by Remark \ref{E} and additivity, both formulae
(\ref{chic}) and (\ref{chismp}) above are satisfied by the
Hodge-Deligne polynomial $E_c(-;u,v)$ defined by means of compactly
supported cohomology. In other words, the polynomial $E_c(-;u,v)$
satisfies the stratified multiplicative property.
\end{remark}

\subsection{$\chi_y$-genera and mixed Hodge modules.}\label{calc} In this
section, by using Saito's theory of mixed Hodge modules we derive
some easy additivity properties of $\chi_y$-genera of complexes of
mixed Hodge structures.  We begin with a consequence of the fact
that mixed Hodge modules over a point are just (graded-polarizable)
mixed Hodge structures:
\begin{lemma}\label{Ho} Let $Z$ be a complex algebraic variety, and $a_Z : Z \to pt$ be the constant map to the point. For any
bounded complex $M^{\bullet}$ of mixed Hodge modules on Z
$$\mathbb{H}^p(Z;M^{\bullet}) :=H^p((a_Z)_* M^{\bullet}) \ \ \ \text{and}
\ \ \  \mathbb{H}_c^p(Z;M^{\bullet}) :=H^p((a_Z)_! M^{\bullet})$$
are mixed Hodge structures. In particular, if $M$ is a mixed Hodge
module on $Z$ whose underlying rational complex is the perverse
sheaf $\mathcal{F}^{\bullet}$, then the hypercohomology group
$\mathbb{H}^p(Z;\mathcal{F}^{\bullet})$ (and resp.
$\mathbb{H}_c^p(Z;\mathcal{F}^{\bullet})$) is the rational vector
space underlying $\mathbb{H}^p(Z;M)$ (and resp.
$\mathbb{H}_c^p(Z;M)$), hence the former is a rational mixed Hodge
structure.
\end{lemma}

As a corollary, we obtain some very useful facts for the
global-to-local study of $\chi_y$-genera. If $Z$ is a complex
algebraic variety, and $i : Y \hookrightarrow Z$ is a closed
immersion, with $j:U \to Z$ the inclusion of the open complement,
then there is a functorial distinguished triangle for $M^{\bullet}
\in D^bMHM(Z)$, lifting the corresponding one from $D^b_c(Z)$ (cf
\cite{Sa1}, p. 321):
\begin{equation}\label{dt} j_!j^*M^{\bullet} \to M^{\bullet} \to
i_*i^*M^{\bullet} \overset{[1]}{\to}\end{equation} In particular, by
taking hypercohomology with compact supports in (\ref{dt}) and
together with Lemma \ref{Ho}, we obtain that the following long
exact sequence is a sequence in the category of mixed Hodge
structures:
\begin{equation}\label{seq}
\cdots \to \mathbb{H}_c^p(U; j^*M^{\bullet}) \to
\mathbb{H}_c^p(Z;M^{\bullet}) \to \mathbb{H}_c^p(Y; i^*M^{\bullet})
\to \cdots
\end{equation}
Therefore, \begin{equation}\label{ad}
\chi_y([\mathbb{H}_c^{\bullet}(Z;
M^{\bullet})])=\chi_y([\mathbb{H}_c^{\bullet}(Y;
M^{\bullet})])+\chi_y([\mathbb{H}_c^{\bullet}(U; M^{\bullet})]).
\end{equation}
As a corollary of (\ref{ad}), by induction on strata we obtain the
following additivity property:
\begin{cor}\label{add2} Let $\mathcal{S}$ be the set of components of strata of
an algebraic Whitney stratification of the complex algebraic variety $Z$.
Then for any $M^{\bullet} \in D^bMHM(Z)$,
\begin{equation}\label{sum}
\chi_y([\mathbb{H}_c^{\bullet}(Z; M^{\bullet})])=\sum_{S \in
\mathcal{S}}\chi_y([\mathbb{H}_c^{\bullet}(S; M^{\bullet})]).
\end{equation}
\end{cor}

\begin{remark}\rm  By taking $M^{\bullet}=\Q_Z^H$
in (\ref{sum}), we obtain the usual additivity of the
$\chi_y^c$-genus (recall that $\chi_y^c$ is defined on the
Grothendieck group of complex varieties). This is a consequence of
the fact that Deligne's and Saito's mixed Hodge structures on
cohomology (with compact support) coincide. The latter assertion can
be seen by construction if the variety can be embedded into a
manifold, but in general it is a very deep result of Saito, see
\cite{Sa5}.
\end{remark}

Each of the terms in the sum of the right-hand side of equation
(\ref{sum}) can be further computed by means of the Leray spectral
sequence for hypercohomology (cf. \S \ref{speq}).

\bigskip

Let $Z$ be an algebraic variety, $M^{\bullet} \in D^bMHM(Z)$, and
fix an algebraic Whitney stratification with respect to which
$\mathcal{F}^{\bullet}:=rat(M^{\bullet}) \in D^b_c(Z)$ has
constructible cohomology. Let $\mathcal{S}$ be the set of components
of strata, and fix $S \in \mathcal{S}$. Then $S$ is a non-singular
complex algebraic variety. Moreover, each cohomology sheaf
$\HC^q(\mathcal{F}^{\bullet})$ on $S$ is a local system underlying a
variation of mixed Hodge structures. We first prove the following
extension of Lemmas \ref{mult} and \ref{cmult} (in fact, the latter
follow from this for some distinguished choices of $M^{\bullet}$):

\begin{prop}\label{mult2}
Assume the local systems $\HC^j(\mathcal{F}^{\bullet})$ are constant
on $S$ for each $j \in \Z$, e.g. $\pi_1(S)=0$. Then
\begin{equation}\label{local}
\chi_y([\mathbb{H}^{\bullet}(S; M^{\bullet})])=\chi_y(S) \cdot
\chi_y([\HC^{\bullet}(\mathcal{F}^{\bullet})_x])
\end{equation}
and
\begin{equation}\label{cslocal}
\chi_y([\mathbb{H}_c^{\bullet}(S; M^{\bullet})])=\chi_y^c(S) \cdot
\chi_y([\HC^{\bullet}(\mathcal{F}^{\bullet})_x]),
\end{equation}
where $[\HC^{\bullet}(\mathcal{F}^{\bullet})_x] \in K_0(mhs)$ is the
complex (with all differentials zero) of mixed Hodge structures
corresponding to stalk cohomologies of $\mathcal{F}^{\bullet}$ at a
point in $S$.
\end{prop}

\begin{proof} This is a consequence of the fact that the spectral
sequences (\ref{leray}) and (\ref{csleray}), applied to the variety
$S$ and to the complex $\mathcal{F}^{\bullet}|_S$, are spectral
sequences of mixed Hodge structures.

By our assumption, the variations of mixed Hodge structures
$\mathcal{H}^q(\mathcal{F}^{\bullet})$ on $S$ are trivial, hence by
(\ref{leray}) and as in the proof of Lemma \ref{mult} there are
mixed Hodge structure isomorphisms
\begin{equation}\label{spt} E_2^{p,q}=H^p(S) \otimes V^q,\end{equation}
where $V^q:= \mathcal{H}^q(\mathcal{F}^{\bullet})_x \cong H^0(S;
\mathcal{H}^q(\mathcal{F}^{\bullet}))$, for any $x \in S$.

Now let $e^{k,l}$ be the Euler-characteristic of Hodge-type $(k,l)$,
i.e., for a bounded complex $K^{\bullet}$ of mixed Hodge structures
$$e^{k,l}([K^{\bullet}])=\sum_i(-1)^ih^{k,l}([K^i]).$$
Since all differentials in (\ref{leray}) are mixed Hodge structure
morphisms, an argument similar to that in Lemma \ref{mult}, yields
that

\begin{eqnarray*} e^{k,l}([\mathbb{H}^{\bullet}(S; M^{\bullet}))]) = \sum_{r+t=k, s+u=l}
e^{r,s}([H^{\bullet}(S)]) \cdot e^{t,u}([V^{\bullet}]).
\end{eqnarray*} Formula (\ref{local}) follows by noting that
$$\chi_y([K^{\bullet}])=\sum_{k,l}e^{k,l}([K^{\bullet}]) \cdot (-y)^k.$$
Formula (\ref{cslocal}) follows similarly, by working instead with
the spectral sequence (\ref{csleray}).

\end{proof}

\begin{remark}\rm It follows from the discussion in \S \ref{def},
that for any $x \in S$, we have:
$$\chi_y([\HC^{\bullet}(\mathcal{F}^{\bullet})_x])=\chi_y([i_x^*\mathcal{F}^{\bullet}]),$$
where $i_x : \{x\} \hookrightarrow S$ is the inclusion of the point.
\end{remark}

Altogether, Corollary \ref{add2} and Proposition \ref{mult2} yield
the following global-to-local formula:
\begin{thm}\label{gl}
Let $\mathcal{S}$ be the set of components of strata of an algebraic
Whitney stratification of the complex algebraic variety $Z$. Assume
that for $M^{\bullet} \in D^bMHM(Z)$ with underlying complex
$\mathcal{F}^{\bullet} \in D_c^b(Z)$, the local systems
$\HC^j(\mathcal{F}^{\bullet})$ are constant on each pure stratum $S
\in \mathcal{S}$ for each $j \in \Z$, e.g. $\pi_1(S)=0$ for all $S
\in \mathcal{S}$. Then
\begin{equation}\label{sum2}
\chi_y([\mathbb{H}_c^{\bullet}(Z; M^{\bullet})])=\sum_{S \in
\mathcal{S}} \chi_y^c(S) \cdot
\chi_y([\HC^{\bullet}(\mathcal{F}^{\bullet})_{x_S}]),
\end{equation}
for some points $x_S \in S$.
\end{thm}

\section{A Hodge-theoretic analogue of the Riemann-Hurwitz formula}\label{curves}
In this section, we apply the above formalism and properties of
$\chi_y$-genera to the study of Hodge-type invariants for a
projective morphism onto a curve. We first recall the definition of
\emph{nearby and vanishing cycle functors}.

\subsection{Vanishing and nearby cycles.}\label{Van}

Let $X$ be a (separated and reduced) complex analytic space of
dimension $n+1$, and $f:X \to \Delta$ be a holomorphic map onto the
unit disc in $\C$, which is smooth over the punctured disc
$\Delta^*$. The total space $X$ is homotopy equivalent to
$X_0=f^{-1}(0)$ by a fiber-preserving retraction $r:X\to X_0$. So
the inclusion $i_t:X_t=f^{-1}(t) \hookrightarrow X$ followed by this
retraction yields the specialization map $r_t :X_t \to X_0$.

Let $i_0:X_0 \hookrightarrow X$ be the inclusion map, and define the
canonical fiber $X_{\infty}$ by
$$X_{\infty}:=X \times_{\Delta^*} \hbar,$$ where $\hbar$ is the
complex upper-half plane (that is, the universal cover of the
punctured disc via the map $z \mapsto \exp(2\pi i z)$). Then
$X_{\infty}$ is homotopy equivalent to any smooth fiber $X_t$ of
$f$. Let $k:X_{\infty} \to X$ be the induced map.
\begin{defn}\rm
The \emph{nearby cycle complex} is defined by
$$\psi_f(\Q_X):=i_0^*Rk_*k^* \Q_X \in D^b_c(X_0).$$
By using a resolution of singularities, it can be shown that
$\psi_f(\Q_X)=R{r_t}_*\Q_{X_t}$, for $t\neq 0$ small enough. The
\emph{vanishing cycle complex} $\phi_f(\Q_X) \in D_c^b(X_0)$ is the
cone on the comparison morphism $\Q_{X_0}=i_0^*\Q_X \to
\psi_f(\Q_X)$, that is, there exists a canonical morphism
$can:\psi_f(\Q_X) \to \phi_f(\Q_X)$ such that
\begin{equation}\label{sp}
i_0^*\Q_X \to \psi_f(\Q_X) \overset{can}{\to} \phi_f(\Q_X)
\overset{[1]}{\to}\end{equation} is a distinguished triangle in
$D^b_c(X_0)$.\end{defn}

In fact, by replacing $\Q_X$ by any complex in $D^b_c(X)$, we obtain
in this way functors $\psi_f, \phi_f:D^b_c(X) \to D^b_c(X_0)$. It is
well-known that if $X$ is locally a complete intersection (e.g., $X$
is smooth), then $\psi_f \Q_X[n]$ and $\phi_f \Q_X[n]$ are perverse
complexes. This is just a particular case of the fact that the
shifted functors $^p \psi_f:=\psi_f[-1]$ and $^p \phi_f:=\phi_f[-1]$
take perverse sheaves on $X$ into perverse sheaves on the central
fiber $X_0$ (cf. \cite{Sc}, Thm. 6.0.2).

\begin{remark}\rm The above construction of the vanishing and nearby
cycles comes up in the following global context (for details, see
\cite{Di2}, \S4.2). Let $X$ be a complex algebraic (resp. analytic)
variety, and $f:X \to \C$ a non-constant regular (resp. analytic)
function. Then for any $t \in \C$, one can construct functors
$$\mathcal{F}^{\bullet} \in D^b_c(X) \mapsto \psi_{f-t}(\mathcal{F}^{\bullet}),
\phi_{f-t}(\mathcal{F}^{\bullet}) \in D^b_c(X_t)$$ where
$X_t=f^{-1}(t)$ is assumed to be a non-empty hypersurface, by simply
repeating the above considerations for the function $f-t$ restricted
to a tube $T(X_t):=f^{-1}(\Delta)$ around the fiber $X_t$ (here
$\Delta$ is a small disc centered at $t$). At least in the algebraic
context, the tube (which is an analytic space) can be chosen so that
$f:T(X_t)\setminus X_t \to \Delta^*$ is a topologically locally
trivial fibration, where $\Delta^*:=\Delta \setminus \{t\}$. In the
analytic context, for this to be true one needs to assume, for
example, that $f$ is proper on the tube $T(X_t)$.
\end{remark}

Of particular importance is the fact that the nearby and vanishing
functors can be defined at the level of Saito's mixed Hodge modules
\cite{Sa0, Sa1}. More precisely, if $f$ is a non-constant regular
(resp. holomorphic) function on the complex algebraic (resp.
analytic) space $X$,  then one has functors $\psi_f, \phi_f :MHM(X)
\to MHM(X_0)$, compatible with the corresponding perverse
cohomological functors on the underlying perverse sheaves by the
functor
$$rat:MHM(X) \to \text{Perv}_{\Q}(X)$$ which assigns to a mixed Hodge module the
underlying $\Q$-perverse sheaf. In other words, $rat \circ \psi_f=
{^p\psi}_f \circ rat$, and similarly for $\phi_f$. As a consequence
of this fact, if $X$ is smooth, for each $x \in X_0$ we have
canonical mixed Hodge structures on the groups
\begin{equation}\label{Mi}
H^j(M_x;\Q)=\mathcal{H}^j(\psi_f \Q_X)_x, \ \
\tilde{H}^j(M_x;\Q)=\mathcal{H}^j(\phi_f \Q_X)_x,\end{equation}
where $M_x$ denotes the Milnor fiber  of $f$ at $x$.

By the identification in (\ref{Mi}), we note that $\text{Supp}
(\phi_f \Q_X) \subset \text{Sing}(X_0)$ (see \cite{Di2}, Example
4.2.6 and Prop. 4.2.8). Then we have the following vanishing result
(e.g., see \cite{Sc}, Example 6.0.13), whose proof we include here
for completeness of the exposition:
\begin{lemma}\label{v} Assume  $X$ is smooth and $f:X \to \Delta$ is proper.
If $s:=\text{dim}_{\C} \text{Sing}(X_0)$, then
$$\mathbb{H}^j(X_0, \phi_f \Q_X)=0, \ \ \text{for} \ j \notin [n-s,
n+s].$$
\end{lemma}
\begin{proof} Since $X$ is a smooth complex manifold of dimension
$n+1$, it follows that $\Q_X[n+1]$ is perverse on $X$. Therefore
${^p\phi}_f\Q_X[n+1]=\phi_f\Q_X[n] \in \text{Perv}(X_0)$. If
$\Sigma:=\text{Sing}(X_0)$, then  $\text{Supp} (\phi_f \Q_X) \subset
\Sigma$ yields that $\phi_f\Q_X[n]|_{\Sigma} \in
\text{Perv}(\Sigma)$, see \cite{Di2}, Cor. 5.2.5. The proof follows
since (noting that $X_0$ and $\Sigma$ are compact)
$$\mathbb{H}^j(X_0, \phi_f \Q_X) \cong \mathbb{H}^{j-n}(X_0, \phi_f
\Q_X[n]) \cong \mathbb{H}^{j-n} (\Sigma, \phi_f \Q_X[n]|_{\Sigma})=0
\ \ \text{for} \ j-n \notin [-s,s],$$ where the vanishing follows
from \cite{Di2}, Proposition 5.2.20.

\end{proof}

\subsection{A Riemann-Hurwitz formula for $\chi_y$-genera.}\label{RH}

Let $f:X \to C$ be a surjective projective  algebraic morphisms from
a smooth $(n+1)$-dimensional complex algebraic variety onto a smooth
algebraic curve. Let $\Sigma(f) \subset C$ be the critical locus of
$f$. Then $f$ is a submersion over $C^*:=C \setminus \Sigma(f)$,
hence locally differentiably trivial (by Ehresmann's fibration
theorem). For a point $c \in C$ we let $X_c$ denote the fiber
$f^{-1}(c)$.

We want to relate the $\chi_y^c$-genera of $X$ and respectively $C$
via the singularities of $f$, and to obtain a stronger version of
Prop. \ref{e} in our setting. The outcome is a Hodge-theoretic
version of a formula of Iversen, or of the Riemann-Hurwitz formula
for the Euler characteristic (e.g., see \cite{Di2}, Cor. 6.2.5,
Remark 6.2.6, or \cite{K}, (III, 32)), see Theorem \ref{Iversen} and
Example \ref{Iv2} below. The proof uses the additivity of the
$\chi_y^c$-genus, together with the study of genera of singular
fibers of $f$ by means of vanishing cycles at a critical value.

\begin{thm}\label{Iversen} Let $f:X \to C$ be a projective algebraic morphism
from a smooth $(n+1)$-dimensional complex algebraic variety onto a
non-singular algebraic curve $C$. Let $\Sigma(f) \subset C$ be the
set of critical values of $f$, and set $C^*=C\setminus \Sigma(f)$.
If the action of $\pi_1(C^*)$ on the cohomology of the generic
fibers $X_t$ of $f$ is trivial (e.g., $\pi_1(C^*)=0$), then
\begin{equation}\label{e1}
\chi_y^c(X)=\chi_y^c(C) \cdot \chi_y^c(X_t) - \sum_{c \in
\Sigma(f)}\chi_y([\mathbb{H}^{\bullet}(X_c;\phi_{f-c}\Q_X)])
\end{equation}
\end{thm}

\begin{proof} Under our assumptions, the fibers of $f$ are complex
projective varieties, and fibers over points in $C^*$ are smooth. By
additivity, we can write:
$$\chi_y^c(X)=\chi_y^c(X^*)+\sum_{c \in \Sigma(f)} \chi_y(X_{c}),$$
where $X^*:=f^{-1}(C^*)$. Then by Lemma \ref{mult}, we have that
$\chi_y^c(X^*)=\chi_y^c(C^*)\chi_y(X_t)$, where $X_t$ is the smooth
(generic) fiber of $f$.

Now let $c \in \Sigma(f)$ be a critical value of $f$ and restrict
the morphism to a tube $T(X_c):=f^{-1}(\Delta_c)$ around the
singular fiber $X_c$, where $\Delta_c$ denotes a small disc in $C$
centered at $c$. By our assumptions, $f :T(X_c) \to \Delta_c$ is a
proper holomorphic function, smooth over $\Delta_c^*$, with fibers
complex projective varieties, that is, a one-parameter degeneration
of complex projective varieties. Then there is a long exact sequence
of mixed Hodge structures (e.g., see \cite{NA}, Thm. 1.1):
\begin{equation}\label{NA}
\cdots \to H^j(X_c;\Q) \to \mathbb{H}^j(X_c;\psi_{f-c}\Q_X) \to
\mathbb{H}^j(X_c; \phi_{f-c} \Q_X) \to \cdots,
\end{equation}
where $\HB^j(X_c;\psi_{f-c}\Q_X)$ carries the \emph{limit mixed Hodge
structure} defined on the cohomology of the canonical fiber (usually
denoted $X_{\infty}$) of the one-parameter degeneration $f :T(X_c)
\to \Delta_c$ (e.g., see \cite{PS}, \S 11.2). However, a consequence
of the definition of the limit mixed Hodge structure is that (cf.
\cite{PS}, Cor. 11.25)
$$\text{dim}_{\C}F^p H^j(X_{\infty};\C)=\text{dim}_{\C}F^p
H^j(X_t;\C),$$ where $X_t$ is the generic fiber of the family (and
of $f$). Therefore,
$$\chi_y(X_{\infty}):=\chi_y([\HB^{\bullet}(X_c;\psi_{f-c}\Q_X)])=\chi_y(X_t).$$
With this observation, from (\ref{NA}) we obtain that for a critical
value $c$ of $f$ the following holds:
\begin{eqnarray*}\label{e3}
\chi_y(X_c) &=&
\chi_y([\mathbb{H}^{\bullet}(X_c;\psi_{f-c}\Q_X)])-\chi_y([\mathbb{H}^{\bullet}(X_c;\phi_{f-c}\Q_X)])
\\ &=&
\chi_y(X_t)-\chi_y([\mathbb{H}^{\bullet}(X_c;\phi_{f-c}\Q_X)]).
\end{eqnarray*}
By additivity and Lemma \ref{mult}, this yields (\ref{e1}).

\end{proof}

\begin{remark}\rm The key point in the proof of the above theorem was to observe that in a
one-parameter degeneration of complex projective manifolds the
$\chi_y$-genus of the canonical fiber coincides with the
$\chi_y$-genus of the generic fiber of the family. Note that this
fact is not true for the corresponding $E$-polynomials, since, while
the Hodge structure on the cohomology of the generic fiber is pure,
the limit mixed Hodge structure on the cohomology of the canonical
fiber carries the monodromy weight filtration.
\end{remark}

\begin{example}\label{Iv2}\rm If $X$ is smooth and $f$ has only isolated singularities, then
\begin{equation}\label{e2}
\chi_y^c(X)=\chi_y^c(C)\chi_y^c(X_t) + (-1)^{n+1} \sum_{x \in
\text{Sing}(f)}\chi_y([\tilde{H}^n(M_x;\Q)]),
\end{equation}
where $M_x$ is the Milnor fiber of $f$ at $x$. Indeed, if $f$ has
only isolated singular points, then each critical fiber $X_c$ has
only isolated singularities and the corresponding vanishing cycles
$\phi_{f-c}\Q_X$ are supported only at these points. Then (\ref{e2})
follows from the identification in (\ref{Mi}) and the vanishing of
Lemma \ref{v}.
\end{example}

\begin{remark}\rm
In the special case of the Euler characteristic $\chi=\chi_{-1}$,
the formula in Theorem \ref{Iversen} and in the example above holds
for any proper analytic morphism onto a curve, without any
assumption on the monodromy (see \cite{Di2}, Cor. 6.2.5). This
follows from the multiplicativity of the Euler characteristic $\chi$
under fibrations, the additivity of compactly supported Euler
characteristic $\chi_c$, and from the fact that $\chi=\chi_c$ (cf.
\cite{F}, pp. 141-142).
\end{remark}

We will extend the formula of Example \ref{Iv2} to the case of
general singularities.  By (\ref{e1}), it suffices to restrict $f$
over a small disc $\Delta_c$ centered at a critical value $c \in
\Sigma(f)$ and to study the polynomial
$\chi_y([\mathbb{H}^{\bullet}(X_c;\phi_{f-c}\Q_X)])$. Recall that
the fibers of $f$ are complex projective algebraic varieties, which
are smooth over points in $\Delta_c \setminus \{c\}$.

Fix  an algebraic Whitney stratification of $X_c$ with respect to
which $\phi_{f-c}\Q_X$ is constructible. For each $q \in \Z$ and
each pure stratum $S \subset \text{Sing}(X_c)$,
$\mathcal{H}^q(\phi_{f-c}\Q_X)$ is a local coefficient system on $S$
with stalk $\tilde{H}^q(M_S;\Q)$, where $M_S$ is the local Milnor
fibre at a point in $S$. Then, according to Theorem \ref{gl}, for
$M^{\bullet}=\phi_f\Q_X[n]$ and by assuming trivial monodromy along
all strata $S \subset \text{Sing}(X_c)$, we obtain
$$\chi_y([\mathbb{H}^{\bullet}(X_c;\phi_{f-c}\Q_X)])=\sum_{S
 \subset \text{Sing}(X_c)} \chi_y^c(S) \cdot
\chi_y([\tilde{H}^{\bullet}(M_S;\Q)]),$$ where $M_S$ is the local
Milnor fibre of $f$ at a point in $S$.

All these facts yield the following general result:

\begin{cor}
Let $f:X \to C$ be a projective algebraic morphism from a smooth
$(n+1)$-dimensional complex algebraic variety onto a non-singular
algebraic curve $C$. Let $\Sigma(f) \subset C$ be the set of
critical values of $f$, and set $C^*=C\setminus \Sigma(f)$. Assume
each special fiber $X_c$ has an algebraic stratification with
respect to which the corresponding vanishing cycle complex is
constructible, and moreover the monodromy along each pure stratum is
trivial. If the action of $\pi_1(C^*)$ on the cohomology of the
generic fibers $X_t$ of $f$ is trivial, then
\begin{equation}\label{e5}
\chi_y^c(X)=\chi_y^c(C) \cdot \chi_y^c(X_t) - \sum_{c \in \Sigma(f)}
\sum_{S \subset \text{Sing}(X_c)} \chi_y^c(S) \cdot
\chi_y([\tilde{H}^{\bullet}(M_S;\Q)])
\end{equation}

\end{cor}


\section{The presence of monodromy. Atiyah-Meyer formulae for the
$\chi_y$-genus.}\label{nm}

In this section we prove Hodge-theoretic analogues (in the category
of complex algebraic varieties) of the Atiyah formula for the
signature of a fibre bundles in the presence of monodromy \cite{At},
and of Meyer's twisted signature formula \cite{Me}. We first state
and prove our formulae in the context of smooth projective
varieties, then point out several interesting extensions to more
general situations.

\subsection{Hirzebruch classes of complex projective manifolds and the Hirzebruch-Riemann-Roch theorem.}
Recall that if $X$ is a smooth complex projective variety, its
Hirzebruch class $\widetilde{T}_y^*(T_X)$ corresponds to the
(un-normalized) power series $$\widetilde{Q}_y(\alpha):=\frac{\alpha
(1+ye^{-\alpha})}{1-e^{- \alpha}}\in \Q[y][[\alpha]], \ \ \
\widetilde{Q}_y(0)=1+y.$$ In fact,
$$\widetilde{T}_y^*(T_X):=td^*(X) \cup ch^*(\lambda_y(T^*_X)),$$ where
$td^*(X)$ is the total Todd class of $X$, $ch^*$ is the Chern
character, and $\lambda_y(T^*_X):=\sum_p \Lambda^pT^*_X \cdot y^p$
is the total $\lambda$-class of (the cotangent bundle of) $X$.
Hirzebruch's class appears in the \emph{generalized
Hirzebruch-Riemann-Roch theorem} (cf. \cite{H}, \S 21.3, but for the
version needed here see also \cite{Y,BSY}), which asserts that if
$E$ is a holomorphic vector bundle on $X$ then the
$\chi_y$-characteristic of $E$, which is defined by
$$\chi_y(X,E):=\sum_{p \geq 0} \chi(X, E \otimes \Lambda^pT^*_X)
\cdot y^p=\sum_{p \geq 0} \left( \sum_{i \geq 0} (-1)^i \text{dim}
H^i(X,\Omega(E) \otimes \Lambda^pT^*_X) \right) \cdot y^p,$$ with
$T^*_X$ the cotangent bundle of $X$ and $\Omega(E)$ the coherent
sheaf of germs of sections of $E$\footnote{ For $X$ smooth and
projective, $\chi_y(X,\mathcal{O}_X)$ agrees with the
Hodge-theoretic $\chi_y$-genus defined in the first part of this
paper. Indeed, by Deligne's theory, one has the following equality
for the Hodge numbers: $h^{p,q}=\text{dim}_{\C} H^q(X,
\Lambda^pT^*_X)$.}, can in fact be expressed in terms of the Chern
classes of $E$ and the tangent bundle of $X$, or more precisely
\begin{equation}\label{gHRR}
\chi_y(X,E)=\int_X \left( ch^*(E) \cup \widetilde{T}_y^*(T_X)
\right) \cap [X].
\end{equation}
In particular, if $E=\mathcal{O}_X$ we have that $$\chi_y(X)=\int_X
\widetilde{T}_y^*(T_X) \cap [X].$$ Also note that the value $y=0$ in
(\ref{gHRR}) yields the classical Hirzebruch-Riemann-Roch theorem
(in short, HRR) for the holomorphic Euler characteristic of $E$,
that is (cf. \cite{H}), $$\chi(X,E)=\int_X \left( ch^*(E) \cup
td^*(X) \right) \cap [X].$$

\subsection{$\chi_y$-genera of smooth projective families.}\label{smfam}
Let $f:E \to B$ be a smooth proper map of smooth complex projective
varieties. By Ehresmann's theorem, $f$ is a differentiable
fibration. The cohomology groups $H^k(E_b)$ of fibers fit into a
local system $R^kf_*\Q_E$. By a result of Griffiths, this local
system underlies a weight $k$ (geometric) variation of Hodge
structures on $B$ such that the Hodge structure at $b \in B$ is just
the Hodge structure we have on $H^k(E_b)$. Let $$\HC_k:=R^kf_*\Q_E
\otimes_{\Q} \mathcal{O}_B.$$ This is a holomorphic bundle with a
flat connection $\bigtriangledown: \HC_k \to \HC_k
\otimes_{\mathcal{O}_B} \Omega^1_B$, and it admits a finite
decreasing filtration $\{\mathcal{F}^p_k\}_p$ by holomorphic
sub-bundles satisfying Griffiths' transversality condition
$\bigtriangledown (\mathcal{F}_k^p) \subset \mathcal{F}_k^{p-1}
\otimes \Omega^1_B$. Set $\HC^{p,k-p}:=Gr^p_{\mathcal{F}}\HC_k$.
This is a holomorphic bundle, not necessarily flat, of rank
$h^{p,k-p}_b:=\text{dim} H^{p,k-p}(E_b;\C)$. The numbers
$h^{p,k-p}_b$ for $b \in B$, remain constant in the family since
they depend upper-semicontinuously on $b$.\newline

One of the main results of this section is the following:
\begin{thm}\label{AM} Let $f:E \to B$ be a smooth proper map of smooth complex projective
varieties. Then the Hirzebruch $\chi_y$-genus of $E$ can be computed
by the following formula:
\begin{equation}\label{twisted}
\chi_y(E)=\int_B \left( ch^*\left( \chi_y (f) \right) \cup
\widetilde{T}^*_y(T_B) \right) \cap [B],
\end{equation}
where $\chi_y(f):=\sum_{p,q \geq 0} (-1)^q \HC^{p,q} \cdot y^p \in
K(B)[y]$ is the $K$-theory $\chi_y$-genus of $f$.
\end{thm}

Before proving the theorem, we make few remarks and state some
immediate consequences.

\begin{remark}\rm \label{high}
\medskip

\noindent$(1)$ Formula (\ref{twisted}) shows the deviation from
multiplicativity of the $\chi_y$-genus of fiber bundles in the
presence of monodromy. The right-hand side of (\ref{twisted}) is a
sum of polynomials, one of the summands being $\chi_y(B) \cdot
\chi_y(F)$. Indeed, the zero-dimensional piece of $ch^*(\chi_y(f))$
is $\chi_y(F)$.

\medskip

\noindent$(2)$ Formula (\ref{twisted}) is the Hodge-theoretic
analogue of Atiyah's signature formula (cf. \cite{At}, (4.3)) in the
complex algebraic setting. Indeed, if $y=1$, then by \cite{Y},
Remark 3, and in the notation of \cite{At} \S 4,
$$\widetilde{T}_1^*(T_B)=\prod_{i=1}^{dim B}
\frac{\alpha_i}{\tanh(\frac{1}{2}\alpha_i)}=:\widetilde{\mathcal{L}}(B),$$
where $\alpha_i$ are the Chern numbers of the tangent bundle of $B$.
Moreover, it is known that $\chi_1(E)=\sigma(E)$ is the usual
signature (cf. \cite{H}), and in a similar fashion one can show that
$(-1)^q \HC^{p,q}$ is the $K$-theory signature $\text{Sign}(f)$ from
\cite{At}. In other words, the value at $y=1$ of (\ref{twisted})
yields $$\sigma(E)=\int_B \left( ch^*(\text{Sign}(f)) \cup
\widetilde{\mathcal{L}}(B) \right) \cap [B].$$

\medskip

\noindent$(3)$ In \cite{At}, Atiyah pointed out that the
non-multiplicativity examples for the signature of holomorphic
fibrations $Z \to C$ having $\text{dim}_{\C} Z=2$, $\text{dim}_{\C}
C=1$, with non-trivial monodromy action on the cohomology of the
fiber $G$, also show the non-multiplicativity of the Todd genus.
Examples of non-multiplicativity in higher dimensions can be
obtained as follows. Let $D \to C$ be an arbitrary holomorphic fiber
bundle with fiber $F$ and having a trivial monodromy. Then the Todd
(and hence $\chi_y$) genus is non-multiplicative for the fibration
$Z \times_C D \to D$, since $Z \times_C D $ also fibers over $Z$
with fiber $F$ and trivial monodromy, and
\begin{equation}
Td(Z \times_C D)=Td(F)Td(Z)\ne Td(F)Td(C)Td(G)=Td(D)Td(G)
\end{equation}
More such examples can be obtained via standard constructions, e.g.
fiber or direct products of Atiyah examples, or higher dimensional
examples  as above.

\medskip

\noindent$(4)$ Theorem \ref{AM} can be extended so that we allow $E$
and $F$ to be singular. We can also discard the compactness
assumption on the base $B$, but in this case we need to allow
contributions ``at infinity" in our formula; see the discussion at
the end of this section for precise formulations of these general
results.
\end{remark}

An immediate corollary of Theorem \ref{AM} is the following:
\begin{cor}\label{triv} Under the assumptions of the above theorem, if moreover
$R^kf_*\Q_E$ is a local system of Hodge structures for each $k$,
i.e. the monodromy action of $\pi_1(B)$ on the cohomology of the
fiber preserves the Hodge filtration, then
\begin{equation}\label{multl}
\chi_y(E)=\chi_y(B) \cdot \chi_y(F),
\end{equation}
where $F$ is the typical fiber of the family.
\end{cor}
\begin{proof}
Indeed, if $R^kf_*\Q_E$ is a local system of Hodge structures, the
Griffiths transversality condition for the flat connection
$\bigtriangledown: \HC_k \to \HC_k \otimes_{\mathcal{O}_B}
\Omega^1_B$ reduces to $\bigtriangledown (\mathcal{F}_k^p) \subset
\mathcal{F}_k^{p} \otimes \Omega^1_B$. It follows that all bundles
$\mathcal{F}_k^p$, whence the bundles
$Gr^p_{\mathcal{F}}\HC_k=\HC^{p,k-p}$, are flat. Since the rational
Chern classes in positive degrees of flat bundles are trivial (cf.
\cite{KT}), we obtain
$$ch^*\left(\sum_{p,q} (-1)^q \HC^{p,q} y^p \right) = \sum_{p,q} (-1)^q \text{rank} (\HC^{p,q}) y^p=\chi_y(F).$$
The result follows.

\end{proof}

\begin{remark}\label{fail}\rm
\medskip
\noindent$(1)$ The above corollary is false in the non-compact case,
e.g., for the Milnor fibration of the cuspidal cupid that was
already considered in \S \ref{lem} (in general, the monodromy of a
weighted homogeneous hypersurface singularity is an algebraic
morphism, thus induces a morphism of mixed Hodge structures in
cohomology).

\medskip

\noindent$(2)$ An example of fibration with the action as in
Corollary \ref{triv} can be given as follows. Let $G$ be a finite
group of biholomorphic maps acting freely on the smooth projective
varieties $E$ and $F$, so that the action on the cohomology of $F$
is non-trivial. If we let $G$ act diagonally on $E \times F$, then
the fibration:
\begin{equation}
(E \times F)/G \to E/G
\end{equation}
has $F$ as its fiber, the monodromy action coincides with the action
of $G$  on the cohomology of $F$, and it preserves the Hodge
filtration on $H^*(F)$ since the monodromy transformations are
biholomorphic.
\end{remark}

\begin{remark}\label{highg}\rm\ \emph{Higher $\chi_y$-genera.} If $X$ is a smooth projective
variety, $\pi:=\pi_1(X)$, and $\alpha \in H^*(B\pi;\Q)$, we define
higher $\chi_y$-genera of $X$ by the formula
\begin{equation}\label{hi}
\chi_y^{[\alpha]}(X):=\int_X \left( u^*(\alpha) \cup
\widetilde{T}_y^*(T_X) \right) \cap [X],
\end{equation}
where $u:X \to B\pi$ is the classifying map of the universal cover
of $X$.

Under the assumptions of Corollary \ref{triv}, formula
(\ref{twisted}) can be rephrased in terms of higher $\chi_y$-genera
as follows. From the classifying space description of each of the
bundles $\HC^{p,q}$, it is clear that $ch^*\left( \HC^{p,q} \right)$
is induced from an universal characteristic class $ch^*(H^{p,q}) \in
H^*(BGL(r;\C);\Q)$, where $r=\text{rank} \HC^{p,q}$. Moreover, the
assumption that the action of $\pi_1(B)$ preserves the Hodge
filtration, hence the $(p,q)$-type, yields that the classifying map
$B \to BGL(r;\C)$ for $\HC^{p,q}$ factors (up to homotopy) as  $B
\overset{u}{\to} B\pi \overset{v}{\to} BGL(r;\C)$, where $u$ is the
classifying map of the universal cover of $B$, and $v$ is induced by
the monodromy action. It follows that $ch^*\left( \HC^{p,q}
\right)=u^*(v^* ch^*(H^{p,q}))$.  Set $\alpha^{p,q}_{\pi}:=v^*
ch^*(H^{p,q})$. Then the Chern character of the $K$-theory
$\chi_y$-genus of $f$ can be written as
$$ch^*(\chi_y(f))=u^* \left( \sum_{p,q}(-1)^q \alpha^{p,q}_{\pi} \cdot y^p \right)$$
Then in the notation of (\ref{hi}) and under the assumptions of
Corollary \ref{triv}, formula (\ref{twisted}) asserts that
$\chi_y(E)$ can be written as a Hodge polynomial in higher genera of
$B$, namely \begin{equation}\label{high}\chi_y(E)=\sum_{p,q} (-1)^q
\chi_y^{[\alpha^{p,q}_{\pi}]}(B) \cdot y^p.\end{equation} However,
since we assumed that each $\HC^{p,q}$ is a flat bundle, we have
that $\chi_y^{[\alpha^{p,q}_{\pi}]}(B)=h^{p,q}(F) \cdot \chi_y(B)$,
and the equation (\ref{high}) yields the multiplicativity of
Corollary \ref{triv}.

\end{remark}

We now return to the proof of Theorem \ref{AM}.
\begin{proof} Recall from \S \ref{speq} that the Leray spectral sequence of the map
$f$, that is,
\begin{equation} E_2^{p,q}=H^p(B, R^qf_*\Q_E)
\Longrightarrow H^{p+q}(E),\end{equation} is a spectral sequence of
mixed Hodge structures. In fact all mixed Hodge structures involved
in the case of a smooth projective family are pure (and by a result
of Deligne, the Leray spetral sequence of $f$ degenerates at $E_2$).
Indeed, the local systems $R^qf_*\Q_E$ underlie geometric variations
of pure Hodge structures (thus admissible in the sense of
Steenbrink-Zucker-Kashiwara; see \cite{PS}, Theorem 14.49 and the
references therein), and it is known (e.g., from Saito's work
\cite{Sa1}, see also \cite{PS}, Theorem 14.50) that the cohomology
of a smooth projective variety with coefficients in such a variation
admits a pure Hodge structure.

By definition, we have
$$\chi_y(E)=\sum_{i,p}(-1)^i \text{dim} Gr^p_F H^i(E; \C) \cdot
(-y)^p=\sum_p \chi^p(E) \cdot (-y)^p,$$ where we let $ \chi^p(E)
:=\sum_i (-1)^i \text{dim} Gr^p_F H^i(E; \C)$  be the Euler
characteristic associated to the exact functor $Gr^p_F$. Since the
differentials in the Leray spectral sequence of $f$ are morphisms of
(mixed) Hodge structures, thus strict with respect to the Hodge
filtrations, it follows that {\allowdisplaybreaks
\begin{eqnarray*} \chi^p(E) &=& \sum_{k,l}
(-1)^{k+l}\text{dim}Gr^p_FH^k(B, R^lf_*\C_E) \\
&=& \sum_{l} (-1)^l \left( \sum_k (-1)^k \text{dim}Gr^p_FH^k(B,R^lf_*\C_E) \right)\\
&=& \sum_{l} (-1)^l \chi^p(B,R^lf_*\Q_E). \end{eqnarray*} }
Therefore, {\allowdisplaybreaks
\begin{eqnarray*}\label{tot} \chi_y(E) &=& \sum_{p} \left( \sum_{l} (-1)^l \chi^p(B,R^lf_*\Q_E) \right) \cdot
(-y)^p\\
&=& \sum_{l} (-1)^l \left( \sum_p \chi^p(B,R^lf_*\Q_E) \cdot
(-y)^p \right)\\
&=& \sum_{l} (-1)^l \chi_y(B,R^lf_*\Q_E). \end{eqnarray*} }

So, we reduced the problem to the following setting, which is a
Hodge-theoretic analogue of the situation considered by Meyer
\cite{Me}: $B$ is a smooth projective variety, $\VB^l:=R^lf_*\Q_E$
is a geometric variation of pure Hodge structures of weight $l$ on
$B$ (in fact, for what follows, one may replace $\VB^l$ by a weight
$l$ polarized variation of Hodge structures on $B$, or more
generally, by an admissible variation of mixed Hodge structures),
and we consider the $\chi_y$-genus of $B$ twisted by $\VB^l$, that
is, $\chi_y(B,\VB^l)$, which encodes the Hodge numbers of the
(mixed, if $\VB^l$ is replaced by an admissible variation) Hodge
structures on the cohomology groups $H^k(B; \VB^l)$, $k \in \Z_{\geq
0}$.

If we let, as before, $\HC_l:=\VB^l \otimes_{\Q} \mathcal{O}_B$ be
the flat bundle associated to $\VB^l$, then we have an isomorphism
$$H^k(B; \VB^l \otimes \C) \cong \HB^k(B; \Omega_B^{\bullet}
\otimes_{\mathcal{O}_B} \HC_l),$$ and the Hodge filtration on
$H^k(B; \VB^l \otimes \C)$ is induced by the filtration
$F^{\bullet}$ on the de Rham complex that is defined by Griffiths'
transversality
$$F^p(\Omega_B^{\bullet}
\otimes_{\mathcal{O}_B} \HC_l):= \left[ \mathcal{F}_l^{p} \HC_l
\overset{\bigtriangledown}{\to} \Omega_B^1 \otimes
\mathcal{F}_l^{p-1} \HC_l \overset{\bigtriangledown}{\to} \cdots
\overset{\bigtriangledown}{\to} \Omega_B^i \otimes
\mathcal{F}_l^{p-i} \HC_l \overset{\bigtriangledown}{\to} \cdots
\right]$$ The associated graded is the complex $$Gr^p_F
(\Omega_B^{\bullet} \otimes_{\mathcal{O}_B} \HC_l)=\left(
\Omega_B^{\bullet} \otimes _{\mathcal{O}_B}
Gr^{p-{\bullet}}_{\mathcal{F}} \HC_l, \ Gr_F \bigtriangledown
\right)$$ with the induced differential.

Then {\allowdisplaybreaks
\begin{eqnarray*}\label{p} \chi^p(B,\VB^l) &=&
\sum_{k} (-1)^k \text{dim } Gr^p_F H^k(B,\VB^l \otimes \C) \\
&=&  \sum_k (-1)^k  \text{dim } Gr^p_F\HB^k(B; \Omega_B^{\bullet}
\otimes_{\mathcal{O}_B} \HC_l)\\
&\overset{(*)}{=}&  \sum_k (-1)^k  \text{dim } \HB^k(B; Gr^p_F
(\Omega_B^{\bullet}
\otimes_{\mathcal{O}_B} \HC_l) )\\
&=& \chi(B,\Omega_B^{\bullet} \otimes _{\mathcal{O}_B}
Gr^{p-{\bullet}}_{\mathcal{F}} \HC_l),
\end{eqnarray*}
} where $(*)$ follows from [\cite{PS}, Theorem 3.18 (iv)] and the
fact proved by M. Saito that $(\VB^l, \Omega_B^{\bullet}
\otimes_{\mathcal{O}_B} \HC_l)$ is a cohomological (mixed, if
$\VB^l$ is replaced by an admissible variation) Hodge complex in the
sense of Deligne (recall $B$ is smooth and compact; see also
\cite{PS}, Theorem 10.9 for the case of a pure polarized variation).

The last term in the above equality can be computed by using the
invariance of the Euler characteristic under spectral sequences. In
general, if $\mathcal{K}^{\bullet}$ is a complex of sheaves on a
topological space $B$, then there is the following spectral sequence
calculating its hypercohomology (e.g., see \cite{Di2}, \S 2.1):
$$E_1^{i,j}=H^j(B,\mathcal{K}^i) \Longrightarrow \HB^{i+j}(B;\mathcal{K}^{\bullet}).$$
Assuming $\chi(B, \mathcal{K}^{\bullet})$ is defined, it can be
computed by $$\chi(B, \mathcal{K}^{\bullet})=\sum_{i,j} (-1)^{i+j}
\text{dim } H^j(B,\mathcal{K}^i)=\sum_i (-1)^i \chi(B;
\mathcal{K}^i)$$

Therefore the twisted $\chi_y$-genus $\chi_y(B,\VB^l)$ can be
computed as follows (where we neglect the cup product symbol or
replace it by  ``$\cdot$" where there is no danger of confusion):
{\allowdisplaybreaks
\begin{eqnarray*}\label{WM} \chi_y(B,\VB^l) &=&
\sum_{p} \chi^p(B,\VB^l) \cdot
(-y)^p\\
&=& \sum_{p} \chi(B,\Omega_B^{\bullet} \otimes
_{\mathcal{O}_B} Gr^{p-{\bullet}}_{\mathcal{F}} \HC_l) \cdot (-y)^p\\
&=& \sum_{i,p} (-1)^i \chi(B,\Omega_B^i \otimes
Gr^{p-i}_{\mathcal{F}} \HC_l) \cdot (-y)^p\\
&\overset{(HRR)}{=}& \sum_{i,p} (-1)^i \left( \int_B
ch^*(Gr^{p-i}_{\mathcal{F}} \HC_l)  ch^*(\Omega_B^i)  td^*(B) \cap
[B] \right) \cdot (-y)^p\\
&=& \int_B  \sum_{i,p} \left( ch^*(Gr^{p-i}_{\mathcal{F}} \HC_l)
\cdot (-y)^{p-i} \right) \cdot \left( ch^*(\Omega_B^i) td^*(B) \cdot
y^i \right) \cap [B] \\
&=& \int_B   \left( \sum_s ch^*(Gr^{s}_{\mathcal{F}} \HC_l) \cdot
(-y)^{s} \right) \cdot \left( td^*(B) \sum_i  ch^*(\Omega_B^i)
\cdot y^i \right) \cap [B] \\
&=& \int_B   \left( \sum_s ch^*(Gr^{s}_{\mathcal{F}} \HC_l) \cdot
(-y)^{s} \right) \cdot   td^*(B) ch^* \left( \lambda_y(T^*_B)
\right)
\cap [B]\\
&=& \int_B   \left( \sum_s ch^*(Gr^{s}_{\mathcal{F}} \HC_l) \cdot
(-y)^{s} \right) \cdot \widetilde{T}^*_y(T_B)  \cap [B].
\end{eqnarray*}
}

Coming back to the computation of $\chi_y(E)$, we obtain that

{\allowdisplaybreaks
\begin{eqnarray*}\label{final} \chi_y(E) &=& \sum_{l} (-1)^l \chi_y(B,R^lf_*\Q_E)\\
&=& \sum_{l} (-1)^l \int_B   \left( \sum_s ch^*(Gr^{s}_{\mathcal{F}}
\HC_l) \cdot (-y)^{s} \right) \cdot \widetilde{T}^*_y(T_B)  \cap
[B]\\
&=& \int_B  \left( \sum_{l,s} (-1)^l ch^*(Gr^{s}_{\mathcal{F}}
\HC_l) \cdot (-y)^{s} \right) \cdot \widetilde{T}^*_y(T_B)  \cap [B]\\
&=& \int_B  \left( \sum_{l,s} (-1)^l ch^*( \HC^{s,l-s}) \cdot
(-y)^{s} \right) \cdot \widetilde{T}^*_y(T_B)  \cap [B]\\
&=& \int_B \left( \sum_{p,q} (-1)^q ch^*( \HC^{p,q}) \cdot y^{p}
\right) \cdot \widetilde{T}^*_y(T_B)  \cap [B].
\end{eqnarray*}
}

\end{proof}

As an important corollary of the proof of Theorem \ref{AM} we obtain
the following Hodge-theoretic analogue of Meyer's signature formula
\cite{Me}:

\begin{cor}\label{Meyer} Let $Z$ be a smooth projective variety and
$\VB$ be a geometric, or polarized (or more generally, an
admissible) variation of (mixed) Hodge structures on $Z$, with
associated flat bundle with ``Hodge" filtration $(\VV,
\mathcal{F}^{\bullet})$. Then the twisted $\chi_y$-genus
$\chi_y(Z,\VB)$ can be computed by the formula
\begin{equation}\label{WM}
\chi_y(Z,\VB)=\int_Z \left( ch^*(Hc_y(\VV)) \cup \widetilde{T}^*_y(T_Z) \right) \cap [Z],
\end{equation}
where $Hc_y(\VV)$ is the $K$-theory Hodge polynomial characteristic of $\VV$ defined by $$Hc_y(\VV)=\sum_p Gr^{p}_{\mathcal{F}} \VV \cdot
(-y)^{p}.$$
\end{cor}

It is now easy to see, with minor changes in the proof of Theorem
\ref{AM}, that we have in fact the following general result:

\begin{thm}\label{gAM} Let $f:E \to B$ be a quasi-projective surjective morphism of complex
algebraic varieties, with $B$ smooth and projective. Assume that the
sheaves $R^sf_*\Z_E$, $s \in \Z$ are locally constant on $B$. Then
the $\chi_y$-genus of $E$ can be computed by the following formula:
\begin{equation}\label{gtwisted}
\chi_y(E)=\int_B \left( ch^*\left( \chi_y (f) \right) \cup
\widetilde{T}^*_y(T_B) \right) \cap [B],
\end{equation}
where $\chi_y(f):=\sum_{i,p} (-1)^i Gr^p_{\mathcal{F}} \HC_i \cdot
(-y)^p \in K(B)[y]$ is the $K$-theory $\chi_y$-genus of $f$.
\end{thm}

\begin{proof} Indeed, by our assumptions of $f$, the local systems
$\LL_s:=R^sf_*\Q_E$, $s \in \Z$, underlie geometric, hence
admissible variations of mixed Hodge structures (see \cite{PS}, Thm.
14.49). Therefore, formula (\ref{WM}) applies to each of these
variations. The rest follows from the Leray spectral sequence of the
map $f$ (cf. (\ref{mapleray})), by noting as in the proof of Thm.
\ref{AM}, that $\chi_y(E)=\sum_s (-1)^s \chi_y(B; \LL_s)$.

\end{proof}

\begin{remark}\rm As stated in
\cite{Me}, Meyer's formula for the signature $\sigma(Z;\mathcal{L})$
of a Poincar\'e local system $\mathcal{L}$ (that is, a local system
equipped with a nondegenerate bilinear pairing $\mathcal{L} \otimes
\mathcal{L} \to \R_Z$)  on a closed oriented smooth manifold $Z$ of
even dimension involves a twisted Chern character and the total
$L$-polynomial of $Z$ (as opposed to Atiyah's formula \cite{At},
where an un-normalized version of the $L$-polynomial is used). More
precisely (\cite{Me}),
$$\sigma(Z;\mathcal{L})=\int_Z \left(
\widetilde{ch^*}([\mathcal{L}]_K) \cup L(Z) \right) \cap [Z],$$
where $[\mathcal{L}]_K$ is the $K$-theory signature of
$\mathcal{L}$, $L(Z)$ is the total Hirzebruch $L$-polynomial of $Z$,
and $\widetilde{ch^*}:=ch^* \circ \psi^2$ is a modified Chern
character obtained by composition with the second Adams operation.
Similarly, following \cite{HBJ}, p.61--62 (see also \cite{SY}, \S
6), we can reformulate our Hodge-theoretic Atiyah-Meyer formulae in
terms of the normalized Hirzebruch classes $T_y^*(T_Z)$
corresponding to the power series
$$Q_y(\alpha):=\widetilde{Q}_y(\alpha (1+y)) \cdot (1+y)^{-1}=\frac{\alpha (1+y)}{1-e^{-\alpha (1+y)}}-\alpha y \in
\Q[y][[\alpha]],$$ by using instead a modified Chern character,
$ch^*_{(1+y)}$, whose value on a complex vector bundle $\xi$ is
$$ch^*_{(1+y)}(\xi)=\sum_{j=1}^{rk \xi} e^{\beta_j (1+y)},$$ for
$\beta_j$ the Chern roots of $\xi$. (In this notation, Meyer's
modified Chern character is simply $ch^*_{(2)}$.) For example, in
the notations of Corollary \ref{Meyer}, formula (\ref{WM}) is
equivalent to
\begin{equation}\label{WMn}
\chi_y(Z,\VB)=\int_Z \left( ch^*_{(1+y)}(Hc_y(\VV)) \cup
{T}^*_y(T_Z) \right) \cap [Z].
\end{equation}

\end{remark}

\bigskip

A similar formula can be obtained for $\chi_y(U;\VB)$, the twisted
$\chi_y$ polynomial associated to the canonical mixed Hodge
structure on $H^*(U;\VB)$, for $U$ any smooth (not necessarily
compact) complex variety and $\VB$ an admissible variation of mixed
Hodge structure on $U$. (The existence of such mixed Hodge
structures follows for example from Saito's theory, see also
\cite{PS}, Thm. 14.50 and the references therein.) In this case, in
order to obtain a cohomological mixed Hodge complex whose Hodge
filtration induces the Hodge filtration on $H^*(U;\VB \otimes \C)$,
we need to use the twisted logarithmic de Rham complex associated to
the Deligne extension of $\VB$ on a good compactification of $U$.
More precisely, let $(\mathcal V, \bigtriangledown)$ be the
corresponding vector bundle on $U$ with its flat connection and
Hodge filtration. Then we can choose a smooth compactification $j:U
\to Z$ such that $D=Z \setminus U$ is a divisor with normal
crossings, and for each half-open interval of length one there is a
unique extension of $(\mathcal V, \bigtriangledown)$ to a vector
bundle $(\bar{\mathcal V}^I, \bar{\bigtriangledown}^I)$ with a
logarithmic connection on $Z$ such that the eigenvalues of the
residues lie in $I$ (cf. \cite{De2}). If we set $\bar {\mathcal
V}:=\bar{\mathcal V}^{[0,1)}$, then the twisted logarithmic de Rham
complex $\Omega_Z^{\bullet}(\text{log}D) \otimes \bar{\mathcal V}$
is quasi-isomorphic (on $Z$) to $Rj_*\VB \otimes \C$, and the
filtration $\mathcal{F}^{\bullet}$ on $\mathcal{V}$ extends to a
filtration $\bar{\mathcal{F}}^{\bullet} \subset \bar{\mathcal V}$
since the variation of mixed Hodge structures was assumed to be
admissible. As before, by Griffiths' transversality, we can filter
the logarithmic twisted de Rham complex,  and Saito proved that this
becomes part of a cohomological mixed Hodge complex that calculates
$H^*(U;\VB)$. By repeating the arguments in the proof of Theorem
\ref{AM}, we obtain the following formula, analogous to (\ref{WM}),
involving contributions at infinity (i.e., forms on $Z$, with
logarithmic poles along $D$):
\begin{equation}\label{nc}
\chi_y(U;\VB)=\int_Z \left(   ch^* (Hc_y(\bar{\VV})) \cup
ch^*(\sum_i \Omega_Z^{i}(\text{log}D) \cdot y^i) \cup td^*(Z)
\right) \cap [Z].
\end{equation}
This explains why under the assumptions of Corollary \ref{triv}, the
multiplicativity of the $\chi_y$-genus fails in the non-compact case
(cf. Remark \ref{fail}).

In view of formula (\ref{nc}), we can obtain an even more general
Atiyah type result for an algebraic map $f$ as in Theorem \ref{gAM}
by dropping the compactness assumption on its target $B$. We leave
the details as an exercise for the interested reader.

\subsection{Higher $\chi_y$-genera and period domains.
}\label{periods}

In this section, we define higher $\chi_y$-genera of variations of
Hodge structures which correspond to cohomology classes of the
quotients of Griffths period domains (cf. \cite{G}). These higher
genera are analogous to the previously considered higher genera
corresponding to the cohomology classes of the fundamental group
(cf. Remark \ref{highg}), and for some types of variations of Hodge
structures coincide with the latter. These classes allow to obtain a
formula for the $\chi_y$-genus of a fibration in terms of
characteristic classes of the base, which yields the
multiplicativity in a variety of cases including the case of trivial
monodromy group.

Let $(B,\VB)$ be a pair where $B$ is a K\"ahler manifold and $\VB$
is a (integer) polarized variation of pure Hodge structures of
weight $k$. If $V$ is the stalk of $\VB$ at a point in $B$, let
$\epsilon=\pm 1$ be the type of the bilinear form $Q$ on $V$ (i.e.
$Q(x,y)=\epsilon Q(y,x)$), and $\eta_V$ be the partition $\text{dim}
V=\sum_{p+q=k} h^{p,q}$ where $h^{p,q}=\text{dim} H^{p,q}$. Let
$D_{\eta_V}$ be the classifying space of pure Hodge structures of
type $(\epsilon,\eta_V)$. This space is a subset in the flag
manifold (consisting of flags in $V$ satisfying the Riemann bilinear
relations) and in particular it is the base of the universal flag
bundle ${\mathcal F}_{\eta_V}:=\cdots \subset \mathcal{F}^p \subset
\cdots$ ($\text{rank} \mathcal{F}^p=\sum_{i \leq p} \text{dim}
H^{i,k-i}$) having the flags as its fiber, and also for the bundles
${\mathcal H}^{p,k-p}$ which are the quotients
$\mathcal{F}^p/\mathcal{F}^{p+1}$.

Let $\bar \Gamma$ be the monodromy group corresponding to $\VB$.
This is a subgroup in the subgroup of $GL(\text{dim} V,{\bf Z})$
consisting of transformations preserving $Q$. The group $\bar
\Gamma$ acts on $D_{\epsilon,\eta_V}$ and some subgroup $\Gamma
\subseteq \bar \Gamma$ of finite index acts freely on
$D_{\epsilon,\eta_V}$. The pair $(B,\VB)$ defines the period map:
$$\pi: B \rightarrow  D_{\epsilon,\eta_V}/\Gamma.$$
The action of the group $\Gamma$ on $V$ and $D_{\epsilon,\eta_V}$
induces the action on the total space of $\mathcal F_{\eta_V}$ so
that the projection ${\mathcal F}_{\eta_V} \rightarrow
D_{\epsilon,\eta_V}$ is $\Gamma$-equivariant. The latter map induces
the locally trivial fibration: ${\mathcal F}_{\eta_V}/\Gamma
\rightarrow D_{\epsilon,\eta_V}/\Gamma$ and, moreover, for any
$h^{p,q} \in \eta_V$ the bundle $\mathcal{H}^{p,q}$ over
$D_{\epsilon,\eta_V}$ descends to a bundle over the quotient
$D_{\epsilon,\eta_V}/\Gamma$.

\begin{defn}\label{highgd} Let $\alpha \in H^*(D_{\epsilon,\eta_V}/\Gamma)$. The higher genus
$\chi_y^{[\alpha]}$ is given by:
$$\chi_y^{[\alpha]}=\pi^*(\alpha) \cup \widetilde{T}_y(B)[B].$$
\end{defn}

Among variations of Hodge structures one can single out those for
which, in the case $\epsilon=-1$ there are at most two non-vanishing
Hodge numbers, and for $\epsilon=+1$ if $p \ne q$ then all
$h^{p,q}=0$ except for at most two of them for which one has
$h^{p,q}=1$. In this case, the period domain is simply-connected
since it is the Siegel upper-half plane for $\epsilon=-1$ and
$SO(2,h^{p,p})/U(1) \times SO(h^{p,p})$, i.e. the quotient by the
maximal compact subgroup, for $\epsilon=+1$ (cf. \cite{CMSP},
p.145). In this case the period map factors as $B \rightarrow
B\pi_1(B) \rightarrow D_{\epsilon,\eta_V}/\Gamma =B(\Gamma)$ (the
latter is the classifying space of $\Gamma$), and the
$\chi_y^{[\alpha]}$ coincides with the higher-$\chi_y$ genus
considered in \cite{BL} (see also Remark \ref{highg}). We shall
refer to such variations as topological variations of Hodge
structures.

Next we shall consider a geometric variation of Hodge structures
arising from a smooth proper map of smooth projective varieties $f:E
\rightarrow B$ (or of compact K\"ahler manifolds). The K\"ahler
class $\omega_E$ of $E$ induces the K\"ahler class $\omega_F$ on
each fiber $F$, which is left invariant by the monodromy of this
fibration. In particular the Hodge form on $H^k(F)$ (given by
$Q(\alpha, \beta)= (\omega_F^{dim F-2k} \cup \alpha \cup \beta)[F]$)
is a monodromy invariant. The primitive cohomology of $F$ yields a
polarized variation of pure Hodge structures, and the fundamental
group $\pi_1(B)$ acts via the monodromy representation on $\oplus_k
H^k_{prim}(F)$ and also on $\prod_k D_{(-1)^k,\eta_k}$, where
$\eta_k$ is the partition of
$\text{dim}H^k(F)_{prim}=\sum_{p+q=k}\text{dim} H^{p,q}(F) \cap
H^k(F)_{prim}$. Let $\bar \Gamma$ be the quotient on $\pi_1(B)$ by
the kernel of this action, and $\Gamma$ be a subgroup of finite
index acting freely.

\begin{thm} Let $f: E \rightarrow B$ be a smooth proper map of smooth projective
varieties, and
$$\pi: B \rightarrow (\prod_k D_{(-1)^k,\eta_k})/\Gamma$$
be the total period map. Let $\pi_k$ be the projection of the target
of $\pi$ on the $k$-th component, and let
$\alpha^{p,q}_{\Gamma}=\pi_k^*(ch^*({\mathcal H}^{p,q}))$ be the
pull back to the quotient of the total period map of the Chern
character of the bundle ${\mathcal H}^{p,q}$. Then $\chi_y(E)$ is
given by the formula (compare with (\ref{high}))
\begin{equation}\label{high2}\chi_y(E)=\sum_{p,q} (-1)^q
\chi_y^{[\alpha^{p,q}_{\Gamma}]}(B) \cdot y^p.\end{equation}
\end{thm}

\noindent The proof follows from the formula (\ref{twisted})
similarly to the way (\ref{high}) was derived from (\ref{twisted}).

\begin{remark}\rm  If $\Gamma=1$ (i.e., the monodromy group $\bar{\Gamma}$ is trivial or
finite) we obtain multiplicativity. More generally, if the $H^{p,q}$
are monodromy invariant then the period map is homotopic to the map
to a point and again one has multiplicativity.\end{remark}

\begin{remark}\rm  Fibrations for which the fibers are curves or K3 surfaces induce
topological variations of Hodge structures and hence the
$\chi_y$-genus of the total space can be expressed in terms of
Novikov's type $\chi_y$-genus as in Remark \ref{highg}. On the other
hand, for fibrations with fibers of higher dimensions one needs the
generalization of the higher $\chi_y$-genus except for very special
cases.\end{remark}

\begin{remark}\rm  The generalization of the higher $\chi_y$-genus given in definition
(\ref{highgd}) has the strong birational invariance property of
higher genera corresponding to twisting by the cohomology classes of
the fundamental group. More precisely, we have the K-equivalence
relation among pairs $(B_1,\VB_1)$ and $(B_2,\VB_2)$ generated by
the elementary K-equivalence $f_1: X \rightarrow B_1, f_2: X
\rightarrow B_2$ such that $f_1^*(\VB_1)=f_2^*(\VB_2)$. The
monodromy groups $\Gamma$ and total period domains $D$ of $\VB_1$
and $\VB_2$ are the same, and for any $\alpha$ in $H^*(D/\Gamma)$
the $\chi_y^{[\alpha]}$-genera of $(B_i,\VB_i)$ coincide. This
follows from the push-forward formulas in the same way as in
\cite{BL}. We shall discuss characterization of such generalized
twisted $\chi_y$ invariants elsewhere.
\end{remark}

\section{Atiyah-Meyer type characteristic class formulae.}\label{CLASS} In this section, we
present characteristic class versions of the above Atiyah-Meyer
formulae for the $\chi_y$-genus. The proofs of these characteristic
class formulae are much more involved, and make use of Saito's
theory of mixed Hodge modules and the construction of the motivic
Hirzebruch classes (cf. \cite{BSY}), which we recall here. For full
details on this construction, the reader is advised to consult
\cite{BSY, CMS}.

Let $Z$ be a complex algebraic variety. Then for any $p \in \Z$ one
has a functor of triangulated categories
$$gr^F_pDR: D^bMHM(Z) \to D^b_{coh}(Z)$$ commuting with proper push-down,
where $D^b_{coh}(Z)$ is the bounded derived category of sheaves of
$\mathcal{O}_Z$-modules with coherent cohomology sheaves. If $\Q_Z^H
\in D^bMHM(Z)$ denotes the constant Hodge module on $Z$, and if $Z$
is smooth and pure dimensional then $gr^F_{-p} DR(\Q_Z^H) \simeq
\Omega^p_Z[-p] \in D^b_{coh}(Z)$. The transformations $gr^F_pDR(M)$
are functors of triangulated categories, so they induce functors on
the level of Grothendieck groups. Thus, if $G_0(Z) \simeq
K_0(D^b_{coh}(Z))$ denotes the Grothendieck group of coherent
sheaves on $Z$, we obtain the following group homomorphism commuting
with proper push-down:
\begin{equation}\label{grF}
gr^F_{-*}DR: K_0(MHM(Z)) \to G_0(Z) \otimes \Z[y, y^{-1}],
\end{equation}
$$[M] \mapsto \sum_p \left( \sum_i (-1)^{i} [\HC^i ( gr^F_{-p} DR(M) )] \right) \cdot (-y)^p.$$

We can now make the following definitions (see \cite{BSY, CMS})
\begin{defn}\label{BSYdefn}\rm
The transformation $\widetilde{MHT_y}$ is defined as the composition
of transformations \footnote{The special case of the transformation
$\widetilde{MHT_y}$ at $y=1$ was previously used by Totaro \cite{To}
for finding numerical invariants of singular varieties, more
precisely Chern numbers that are invariant under small
resolutions.}:
\begin{equation}\label{IT}
\widetilde{MHT_y} :=td_* \circ gr^F_{-*}DR: K_0(MHM(Z)) \to
H_{2*}^{BM}(Z) \otimes \Q[y,y^{-1}],
\end{equation}
where $td_*$ is the Baum-Fulton-MacPherson Todd class transformation
\cite{BFM}, which is linearly extended over $\Z[y,y^{-1}]$. Note
that $\widetilde{MHT_y}$ commutes with proper push-forward.
\end{defn}

\begin{remark}\rm Let $K^0(Z)$ be the Grothendieck group of complex algebraic vector bundles on $Z$. If $Z$ an
algebraic manifold, the canonical map $K^0(Z) \to G_0(Z)$ induced by
taking the sheaf of sections is an isomorphism, and the Todd class
transformation of the classical Grothendieck-Riemann-Roch theorem is
explicitly described by $td_*(\cdot)=ch^*(\cdot)  td^*(T_Z) \cap
[Z]$.
\end{remark}

\begin{defn}\label{Ty}\rm
The \emph{Hirzebruch class} of an $n$-dimensional complex algebraic
variety $Z$ is defined by the formula
\begin{equation}
\widetilde{T_y}_*(Z):=\widetilde{MHT_y}([\Q_Z^H]).
\end{equation}
Similarly, if $Z$ is an $n$-dimensional complex algebraic manifold,
and $\VB$ a polarized variation of Hodge structures on $Z$, we
define \emph{twisted Hirzebruch characteristic classes}
\begin{equation}\label{tHc} \widetilde{T_y}_*(Z;
\VB)=\widetilde{MHT_y}([\VB^H]),\end{equation} where
$\VB^H[n]=((\mathcal{V}, \bigtriangledown),
\mathcal{F}_{-{\bullet}}, \VB[n])$ is the smooth mixed Hodge module
on $Z$ corresponding to $\VB$, with
$\mathcal{F}_{-{\bullet}}:=\mathcal{F}^{{\bullet}}$ (e.g., see
\cite{PS}, Thm. 14.30 and the references therein).
\end{defn}

By [\cite{BSY}, Lemma 3.1 and Theorem 3.1], the following
normalization holds: if $Z$ is smooth and pure dimensional, then
$\widetilde{T_y}_*(Z)=\widetilde{T}^*_y(T_Z) \cap [Z]$, thus
$\widetilde{T_y}_*(Z)$ is an extension  to the singular setting of
(the Poincar\'e dual of) the un-normalized Hirzebruch class.

We can now prove the following Meyer-type formula for the twisted
Hirzebruch characteristic classes

\begin{thm}\label{Meyerc}
Let $Z$ be a complex algebraic manifold of pure dimension $n$, and
$\VB$ a polarized variation of Hodge structures on $Z$ with
associated flat bundle with ``Hodge" filtration $(\VV,
\mathcal{F}^{\bullet})$. Then
\begin{equation}\label{WMc}
\widetilde{T_y}_*(Z;\VB)= \left( ch^*(Hc_y(\VV)) \cup
\widetilde{T_y}^*(T_Z) \right) \cap [Z]= ch^*(Hc_y(\VV)) \cap
\widetilde{T_y}_*(Z),
\end{equation}
where $Hc_y(\VV)$ is the $K$-theory Hodge polynomial characteristic
of $\VV$.
\end{thm}

\begin{proof} Let $\VB^H[n]=((\mathcal{V}, \bigtriangledown), \mathcal{F}_{-{\bullet}},
\VB[n])$ be the smooth mixed Hodge module on $Z$ corresponding to
$\VB$, with $\mathcal{F}_{-{\bullet}}:=\mathcal{F}^{{\bullet}}$ the
increasing filtration on the $\mathcal{D}$-module $\mathcal{V}$. It
follows from Saito's work that there is a filtered quasi-isomorphism
between $(DR(\VB^H),F_{-{\bullet}})$ and the usual filtered de Rham
complex $(\Omega_Z^{\bullet}(\mathcal{V}), F^{\bullet})$ with the
filtration induced by Griffiths' transversality, that is:
$$F^p \Omega_Z^{\bullet}(\mathcal{V}): \left[ \mathcal{F}^p
\overset{\bigtriangledown}{\to} \Omega_Z^1 \otimes \mathcal{F}^{p-1}
\overset{\bigtriangledown}{\to} \cdots
\overset{\bigtriangledown}{\to} \Omega_Z^i \otimes \mathcal{F}^{p-i}
\overset{\bigtriangledown}{\to} \cdots \right].$$ Therefore,
{\allowdisplaybreaks
\begin{eqnarray*} \widetilde{T_y}_*(Z;\VB) &=& td_* \left( \sum_p \left( \sum_i (-1)^{i} [\HC^i ( gr^F_{-p} DR(\VB^H) )] \right) \cdot (-y)^p \right) \\
&=& td_* \left( \sum_p \left( \sum_i (-1)^{i} [\HC^i ( gr_F^{p} \Omega_Z^{\bullet}(\mathcal{V}) )] \right) \cdot (-y)^p \right) \\
&=& td_* \left( \sum_p \left( \sum_i (-1)^{i} [\Omega_Z^{i} \otimes Gr_{\mathcal{F}}^{p-i}\mathcal{V}] \right) \cdot (-y)^p \right) \\
&=& \sum_p  \sum_i (-1)^{i} \ td_* ([\Omega_Z^{i} \otimes Gr_{\mathcal{F}}^{p-i}\mathcal{V}])  \cdot (-y)^p  \\
&\overset{(GRR)}{=}& \sum_p  \sum_i (-1)^{i} \ ch^* (\Omega_Z^{i} \otimes Gr_{\mathcal{F}}^{p-i}\mathcal{V})  \cup td^*(Z) \cap [Z] \cdot (-y)^p  \\
&=& \sum_p  \sum_i \left( ch^* (Gr_{\mathcal{F}}^{p-i}\mathcal{V}) \cdot (-y)^{p-i} \right) \cup \left( ch^*(\Omega^i_Z) \cdot y^i \right) \cup td^*(Z) \cap [Z] \\
&=& \left( \sum_q   ch^* (Gr_{\mathcal{F}}^{q}\mathcal{V}) \cdot (-y)^{q} \right) \cup \left(\sum_i ch^*(\Omega^i_Z) \cdot y^i \right) \cup td^*(Z) \cap [Z] \\
&=& ch^*(Hc_y(\mathcal{V})) \cup ch^*(\lambda_y(T^*_Z)) \cup td^*(Z)
\cap [Z] \\
&=& ch^*(Hc_y(\mathcal{V})) \cup \tilde{T_y}^*(T_Z) \cap [Z],
\end{eqnarray*}
}

\noindent where $(GRR)$ means an application of the classical
Grothendieck-Riemann-Roch theorem.

\end{proof}

\begin{remark}\label{geom}\rm More generally, Saito \cite{Sa1} showed that an admissible variation of
mixed Hodge structures on a smooth variety $Z$ (e.g., a geometric
variation of mixed Hodge structures), with underlying local system
$\LL$, gives rise to a smooth mixed Hodge module on $Z$, call it
$\LL^H[\text{dim} Z]$, with $\LL[\text{dim} Z]$ as its underlying
perverse sheaf. Then one can define as in (\ref{tHc}) twisted
Hirzebruch characteristic classes associated to such admissible
variations. Note that Theorem \ref{Meyerc} remains true in this
greater generality, i.e., we can let $\VB$ be a geometric variation
of (mixed) Hodge structures, or more generally, an admissible
variation of mixed Hodge structures. Then if $Z$ is complete,
Corollary \ref{Meyer} can be obtained from Theorem \ref{Meyerc} by
pushing forward to a point via the constant map $Z \to pt$.
\end{remark}

J\"org Sch\"urmann \cite{Sch} communicated to us that the following
Atiyah-type result can be obtained as a direct application of the
Verdier-Riemann-Roch formula for a smooth proper morphism (cf.
\cite{BSY}, Cor. 3.1(3)), if one makes the identification
$\HC^{p,q}\simeq R^qf_*(\Lambda^pT^*_f)$, with $T^*_f$ the dual of
the tangent bundle $T_f$ to the fibers of $f$ (cf. \cite{PS},
Proposition 10.28). However, the proof we give here is based only on
the definition of the Hirzebruch classes and on Theorem \ref{Meyerc}
in the context of geometric variations of Hodge structures.

\begin{thm}\label{AMc} Let $f:E \to B$ be a smooth projective morphism between complex algebraic manifolds. Then the following holds:
\begin{equation}\label{Htwist}
f_*\widetilde{T_y}_*(E)=  ch^*\left( \chi_y (f) \right) \cap
\widetilde{T_y}_*(B),
\end{equation}
where $\chi_y(f)=\sum_{p,q \geq 0} (-1)^q \HC^{p,q} \cdot y^p$ is
the $K$-theory $\chi_y$-genus of $f$.
\end{thm}

\begin{proof} Since $f$ is proper and the transformation
$\widetilde{MHT_y}$ commutes with proper pushdowns, we first obtain
the following: \begin{equation} f_*\widetilde{T_y}_*(E)=
f_*(\widetilde{MHT_y}([\Q^H_E]))=\widetilde{MHT_y}([f_* \Q^H_E]).
\end{equation}
Now let $\tau_{\leq}$ be the natural truncation on $D^bMHM(B)$ with
associated cohomology $H^{\bullet}$.  Then for any complex
$M^{\bullet} \in D^bMHM(B)$ we have the identification (e.g., see
\cite{Di2}, p. 95-96; \cite{Sc}, Lemma 3.3.1)
\begin{equation}
[M^{\bullet}]=\sum_{i \in \Z} (-1)^i [H^i(M^{\bullet})] \in
K_0(D^bMHM(B)) \cong K_0(MHM(B)).
\end{equation}
In particular, if for any $k \in \Z$ we regard
$H^{i+k}(M^{\bullet})[-k]$ as a complex concentrated in degree $k$,
then
\begin{equation}
\left[ H^{i+k}(M^{\bullet})[-k] \right]= (-1)^k
[H^{i+k}(M^{\bullet})] \in K_0(MHM(B)).
\end{equation}
Therefore, if we let $M^{\bullet}=f_*\Q^H_E$, we obtain that
\begin{equation}\label{int} f_*\widetilde{T_y}_*(E)=\sum_{i \in \Z} (-1)^i \widetilde{MHT_y}([H^i(f_* \Q^H_E)])=\sum_{i \in \Z} (-1)^i
\widetilde{MHT_y}( \left[ H^{i+\text{dim}B}(f_*
\Q^H_E)[-\text{dim}B] \right]).
\end{equation}
Note that $H^i(f_* \Q^H_E) \in MHM(B)$ is the smooth mixed Hodge
module on $B$ whose underlying rational complex is (recall that $B$
is smooth)
\begin{equation} \text{rat} (H^i(f_* \Q^H_E))={^p\HC}^i(Rf_*\Q_E)=(R^{i-\text{dim}B}
f_*\Q_E)[\text{dim} B],
\end{equation}
where for the second equality we refer to \cite{PS}, Example 13.20.
In this case, each of the local systems $\LL_s:=R^{s} f_*\Q_E$
underlies a geometric variation of Hodge structures.

Altogether, (\ref{int}) becomes
\begin{equation} f_*\widetilde{T_y}_*(E)=\sum_{i \in \Z} (-1)^i
\widetilde{T_y}_* (B; \LL_i^H),
\end{equation}
where, by analogy with Definition \ref{Ty},
$\LL_i^H[\text{dim}B]:=H^{i+\text{dim}B}(f_* \Q^H_E)$ is the smooth
mixed Hodge module whose underlying perverse sheaf is
$\LL_i[\text{dim}B]$.

Our formula (\ref{Htwist}) follows now from Theorem \ref{Meyerc} and
Remark \ref{geom}.

\end{proof}

\begin{remark}\rm If $B$ in the above theorem is also complete, then by pushing
(\ref{Htwist}) down to a point via the constant map $B \to pt$, we
get back the result of Theorem \ref{AM}.
\end{remark}

An immediate corollary of Theorem \ref{AMc} is the following
extension of [\cite{CMS}, Cor. 3.12], whose proof imitates that of
Corollary \ref{triv}.
\begin{cor} Under the assumptions of the above theorem, if moreover
$R^kf_*\Q_E$ is a local system of Hodge structures for each $k$,
i.e. the monodromy action of $\pi_1(B)$ on the cohomology of the
fiber preserves the Hodge filtration, then
\begin{equation}
f_*\widetilde{T_y}_*(E)=  \chi_y (F)  \widetilde{T_y}_*(B),
\end{equation}
where $F$ is the fiber of $f$.
\end{cor}

\begin{remark}\rm Note that in the proof of Theorem \ref{AMc} we only used the fact
that $f$ is proper (so that the characteristic class transformation
$\widetilde{MHT_y}$ commutes with $f_*$), and that the local systems
$\LL_s:=R^{s} f_*\Q_E$, $s \in \Z$, underly admissible variation of
mixed Hodge structures on the smooth variety $B$. Therefore Theorem
\ref{AMc} admits the following generalization, where we allow the
generic fiber and the domain of the map to be singular:
\begin{thm}\label{gAMc} Let $f:E \to B$ be a projective surjective morphism of complex
algebraic varieties, with $B$ smooth, such that the sheaves
$R^sf_*\Z_E$, $s \in \Z$ are locally constant on $B$. Then
\begin{equation}\label{gtwistedc}
f_*\widetilde{T_y}_*(E)= ch^*\left( \chi_y (f) \right) \cup
\widetilde{T_y}_*(B),
\end{equation}
where $\chi_y(f):=\sum_{i,p} (-1)^i Gr^p_{\mathcal{F}} \HC_i \cdot
(-y)^p \in K(B)[y]$ is the $K$-theory $\chi_y$-genus of $f$.
\end{thm}
\end{remark}

\providecommand{\bysame}{\leavevmode\hbox
to3em{\hrulefill}\thinspace}

\end{document}